\documentclass[10pt]{extarticle}
\usepackage[toc,page]{appendix}
\usepackage{amsmath,amsthm,amssymb,amsbsy,amsopn,mathtools,comment,bbm,lineno}
\usepackage{graphicx,color}
\usepackage[bookmarks=true,bookmarksnumbered=true,bookmarkstype=toc]{hyperref}
\usepackage{mathptmx,bm}
\usepackage{algorithm}
\usepackage{algorithmicx}
\usepackage{caption}
\usepackage{subcaption}

\usepackage{amsmath}
\usepackage{lineno}

\newtheorem{Def}{Definition}[section]
\newtheorem{theorem}[Def]{Theorem}

\newtheorem{assumption}[Def]{Assumption}

\DeclareMathOperator*{\esssup}{ess\,sup}
\DeclareMathOperator*{\essinf}{ess\,inf}

\usepackage[inline]{showlabels}

\theoremstyle{definition}

\makeatletter
\@addtoreset{equation}{section}
\makeatother

\title{Iterative weak approximation and hard bounds for switching diffusion}
\author{\sc Qinjing Qiu\footnote{Email address: qqiu9438@uni.sydney.edu.au. Postal Address: School of Mathematics and Statistics, University of Sydney, NSW 2006, Australia. TEL: +61(0)2-9114-1264 FAX: +61(0)2-9351-4534} \, and Reiichiro Kawai\footnote{%Corresponding author. 
Email address: raykawai@g.ecc.u-tokyo.ac.jp. Graduate School of Arts and Sciences / Mathematics and Informatics Center, The University of Tokyo, Tokyo 153-8902, Japan, and School of Mathematics and Statistics, The University of Sydney, NSW 2006, Australia. This work was partially supported by JSPS Grants-in-Aid for Scientific Research 20K22301 and 21K03347.}}
	\date{}

\begin{document}
%\linenumbers

\maketitle

\begin{abstract}
\noindent
We establish a novel convergent iteration framework for a weak approximation of general switching diffusion.
The key theoretical basis of the proposed approach is a restriction of the maximum number of switching so as to untangle and compensate a challenging system of weakly coupled partial differential equations to a collection of independent partial differential equations, for which a variety of accurate and efficient numerical methods are available.
Upper and lower bounding functions for the solutions are constructed using the iterative approximate solutions.
We provide a rigorous convergence analysis for the iterative approximate solutions, as well as for the upper and lower bounding functions.
Numerical results are provided to examine our theoretical findings and the effectiveness of the proposed framework.


\end{abstract}

\section{Introduction}
%In this paper, we consider a Markov regime-switching diffusion model, which is a natural extension of the classical diffusion model. 
Owing to the increasing demands on regime-switching diffusion from emerging applications in financial engineering and wireless communications, much attention has been given to switching diffusion processes. A salient feature of such systems is that they include both continuous dynamic and discrete events. 

States of the environment process was referred to in \cite{ZY1} as economic circumstances or political regime-switching.
It is therefore theoretically appealing to modify the classical diffusion process so that both expected return and the variance reflect the changes in the external environmental factors.
The modelling framework that is advocated in this paper achieves this.
The rationale is that in the different modes or regimes, the volatility and return rates can be very different. 
Another example is from a wireless communication networks. 
Consider the performance analysis of an adaptive, linear, multi-user detector in a cellular, direct-sequence, code-division, multiple-access wireless network with changing user activity due to an admission or access controller at the base station. Under certain conditions, an optimization problem for the aforementioned application leads to certain issues, which in turn lead to switching diffusion limit \cite{wireless1}.

%In the model, we assume that both the expected increment and the variance are influenced by an external environment process. This type of process is also known as a Markov-modulated process and was studied in Asmussen \cite{AS2} and Reinhard \cite{RM} two decades ago. 
%For example, 
One early attempt at such hybrid models for financial applications can be traced back to the %\cite{ZY2} 
1960s in which both the appreciation rate and the volatility of a stock depend on a continuous-time Markov chain. The introduction of such models makes it possible to describe stochastic volatility in a relatively simple manner. To illustrate, in the simplest case, a stock market may be considered to have two ``modes'' or ``regimes'', up and down, resulting from the state of the underlying economy, the general mood of investors in the market, and so on.
The primary motivation for the generalization is to enhance the flexibility of the model parameter settings for the classical diffusion process. The examples usually given are weather conditions and epidemic outbreaks, even though seasonality would play a role and can probably not be modelled by a Markov regime-switching model.

A variety of numerical methods extend those that already exist in the single-regime setting, but the extensions are rarely straightforward due to the inter-dependence nature of regime switching systems.
Lattice methods in regime-switching settings have been attempted; first in \cite{Bollen} in a two-regime setting and later in the form of the adaptive mesh model \cite{Gao}, the $2k$-branch lattice for a $k$-regime model \cite{Das}, and the trinomial method \cite{Yang}. Mathematical programming is employed for obtaining hard bounds on regime-switching diffusion in \cite{louis, Kwon}.
Approaches based on penalty methods for pricing American options have also been extended in a non-trivial way to regime-switching frameworks. For instance, in \cite{Labahn, Liu1, Liu2}, exponential time difference schemes are applied in combination with a penalty method to approximate a solution to the complex system of partial differential equations with coupled free boundary conditions that arise in the presence of a regime-switching diffusion.
The Fourier space time-stepping method \cite{Surkov} overcomes some of the problems encountered by lattice methods and finite-difference schemes when it comes to option pricing in models in higher dimensions or models containing jumps. This approach is applicable to a wide class of path-dependent options and options on multiple assets in regime-switching diffusion containing jumps.

A Monte Carlo method is presented in \cite{Scherer} for pricing barrier option in regime-switching diffusion, where the author demonstrates the non-trivial modifications that must be made in order to extend standard simulation-based techniques from the single-regime setting to multiple regimes.
Esscher transform methods are employed in \cite{Siu1} and \cite{Siu2} for pricing European and American options in regime-switching diffusion and jump-diffusion settings.
The Canadization method was extended to regime-switching settings to price American options \cite{Boyarchenko1} and double barrier options \cite{Boyarchenko2}.
We mention the radial basis collocation method \cite{Damircheli}, a mesh-free method for pricing American options under regime-switching jump diffusion models.
A MATLAB toolbox specially designed for the estimation, simulation and forecasting of a general Markov regime-switching diffusion is developed in \cite{MSreg28}.
In all the aforementioned studies, a great deal of special attention must be paid to the switching nature.

%Asymptotic properties of the maximum likelihood estimator in regime-switching models were studied in \cite{HK29}.
%A unified analytical approximation framework for computations under general regime-switching diffusion was introduced in \cite{NSY30}. 
%An estimation method by spectral clustering hidden Markov model is proposed in \cite{ZL31}. 
%A global form stochastic maximum principle for a Markov regime switching mean-field model driven by Brownian motions and Poisson jumps is developed in \cite{ZSX32} and a  perturbation approach is applied to approximate the first-order and second-order estimation in \cite{BM33}. 

In this paper, we establish a novel convergent iteration framework on a generalized Feynman-Kac formula under regime-switching, which are the solutions to a set of weakly coupled partial differential equations with initial and boundary value conditions.
We construct our recursive approximation mechanism by restricting the maximum number of switching occurring before the terminal time in the probabilistic representation in a similar manner to \cite{QIU2022109301, qiu2022recursive}.
The significance of the present approach is that it succeeds to untangle and compensate a challenging set of weakly coupled partial differential equations to a set of independent partial differential equations.
%To the best of our knowledge, 
Such untangling of the inter-dependency is quite remarkable, as undoubtedly a variety of accurate and efficient numerical methods are available for the latter untangled partial differential equations.

The rest of the paper is organized as follows.
We present the necessary background materials and problem formulation regarding regime-switching diffusion in Section \ref{section preliminaries}.
We prove in Sections \ref{section iterative weak approximation} and \ref{section monotonic} that such an approximation can be obtained as the solution to a set of independent partial differential equations with additional heat source (the non-homogeneous terms) and additional discounting rates or, equivalently, the approximation admits a probabilistic representation under a probability measure without regime-switching.
Further, in Section \ref{section hard bounding functions IVP}, we derive hard bounding functions for the true solution under regime-switching based on the approximation at each iteration as well as their convergence to the unknown solution.
After briefly discussing the proposed approach in the context the initial-boundary value problem in Section \ref{section initial boundary value problems}, we examine our theoretical findings through numerical results in Section \ref{section numerical example}.
Finally, Section \ref{section concluding remarks} concludes the present study and highlights future research directions.

%To maintain the flow of the paper, we collect all proofs and technical notes in the Appendix.

\section{Preliminaries}\label{section preliminaries}
We begin with the notation that will be used throughout the paper.
We denote by $\mathbb{R}$ the set of real numbers, and by $\mathbb{N}:=\{1,2,\cdots\}$ the set of natural numbers excluding $0$, with $\mathbb{N}_0:=\mathbb{N}\cup \{0\}$.
We reserve $p$, $n$ and $d$ for fixed positive integers.
We denote by $\{R_t: t\geq 0\}$ a one-dimensional Markov process with right-continuous sample paths taking values in the finite-state space $\mathcal{M} := \{1,\cdots,p\}$.
We define an $n$-dimensional stochastic process $\{X_t: t\geq 0\}$ taking values in a suitable domain $D(\subseteq \mathbb{R}^n)$ by the stochastic differential equation:
\begin{equation}\label{XSDE}
dX_{t}=b(t,X_{t},R_{t})dt+\sigma(t, X_{t},R_{t})dW_{t},\quad (X_0, R_0) \in D\times \mathcal{M},
\end{equation}
where $\{W_t: t \geq 0\}$ is a $d$-dimensional standard Brownian motion and the drift coefficient $b:  [0,\infty) \times D \times \mathcal{M} \to  \mathbb{R}^{n}$ and the diffusion coefficient $\sigma:  [0,\infty) \times D \times \mathcal{M} \to  \mathbb{R}^{n \times d}$ are deterministic functions satisfying the usual conditions \cite{SDE1, SDE2}: 
for each $i \in \mathcal{M}$, $\sigma(\cdot,\cdot,i)$ and $b(\cdot,\cdot,i)$ are at most of linear growth and Lipschitz continuous on every compact subset of  $[0,\infty) \times D$. 
Here, the domain $D(\subseteq \mathbb{R}^n)$ depends on the problem setting under consideration.
For example, we may have $D=(0,\infty)$ when the underlying process $\{X_t:t\geq 0\}$ is a regime-switching geometric Brownian motion starting from a positive level, without an attainable boundary.
For problems with attainable boundaries (Section \ref{section initial boundary value problems}), the domain $D$ is a bounded.
We define the first hitting time of the underlying process on the boundary by
\begin{equation}\label{etadef}
\eta_t:=\inf\{s\geq t: X_s \not\in D\},
\quad t\geq 0;
\end{equation} 
with $\eta_t=+\infty$, if the boundary is never attained.
As usual, we denote by $\partial D$ the boundary of the domain $D$, and by $\overline{D}=D\cup\partial D$ the closure.

In this paper, the transition rate of the regime process $\{R_t: t\geq 0\}$ is Markov and can be state dependent, governed by the so-called generator matrix $Q: D \to \mathbb{R}^{p \times p}$, where we denote its $(i,j)$-entry by $q_{ij}({\bf x}):=(Q({\bf x}))_{i,j}$. 
To ensure that the stochastic differential equation (\ref{XSDE}) admits a unique solution,
 %  described in \cite[Theorem 2.1]{YZ1}
we impose the Assumption \ref{qproperty}, the so-called $q$-property \cite{BYZ,YZ1}, on top of the aforementioned usual conditions on the coefficients $b$ and $\sigma$.

\begin{assumption}[]\label{qproperty}
The matrix $Q$ satisfies the following properties:

\noindent
{\bf (a)} $q_{ij}({\bf x})\geq 0$, for all ${\bf x}\in D$, $i,j \in \mathcal{M}$ and $i\neq j$,

\noindent
{\bf (b)}  $q_{ii}({\bf x}) = - \sum_{j \in\mathcal{M}\setminus\{i\}}q_{ij}({\bf x})$, for all ${\bf x}\in D$ and $i \in \mathcal{M}$,

\noindent
{\bf (c)}  $q_{ij}(\cdot)$ is bounded on $D$, for all $i,j\in\mathcal{M}$, and

\noindent
{\bf (d)}  there exists $\alpha\in(0,1)$ such that $q_{ij}(\cdot)$ is $\alpha$-H\"older continuous on every compact subset of $D$, for all $i,j\in\mathcal{M}$.
\end{assumption}

We write $\mathbb{P}_Q$ for the probability measure under which the transition rate of the regime process $\{R_t:t\geq 0\}$ is given by the generator matrix $Q$.
Throughout, we assume that the paired stochastic process $\{(X_t,R_t):t\geq 0\}$ is defined on the probability space $(\Omega,\mathcal{F},\mathbb{P}_Q)$ endowed with the right-continuous filtration $\mathcal{F}=(\mathcal{F}_t)_{t\geq 0}$.
For $i,j \in \mathcal{M}$, $(t,{\bf x})\in[0,\infty)\times D$ and an infinitesimal $\Delta$,  we formally define
\begin{equation}\label{defPR}
\mathbb{P}_Q^{t,{\bf x},i}\left(R_{t+\Delta} = j \right): =
\mathbb{P}_Q\left(R_{t+\Delta} = j  \big| \,R_t=i, X_t= {\bf x} \right)=
\begin{cases}
q_{ij}( {\bf x})\Delta + o(\Delta), &\text{if } i\neq j, \\
1+q_{ii}( {\bf x})\Delta+o(\Delta),&\text{if } i=j.
\end{cases}
\end{equation}
Hereafter, we denote by $\mathbb{E}_Q^{t,{\bf x},i}$ the expectation taken under the probability measure $\mathbb{P}_Q^{t,{\bf x},i}$.
Next, for $t\geq 0$, we define the time of the first switching occur after time $t$ as:
\[
\tau^{(1)}_t:=\inf\{s>t:R_s\neq R_t\}.
\]
Since all entries of the generator matrix $Q$ are bounded, the probability that the Markov process $\{R_s:s\geq 0\}$ switches more than once within an infinitesimal-sized interval $(t,t+\Delta]$ is as negligibly small as $o(\Delta)$, that is,
\[
\mathbb{P}_Q^{t,{\bf x},i}\left(\tau^{(1)}_t>t+\Delta \right) =
\mathbb{P}_Q^{t,{\bf x},i}\left(R_{t+\Delta} = i \right) = 1+q_{ii}( {\bf x})\Delta+o(\Delta),
\]
due to \eqref{defPR}. 
For $m\in\mathbb{N}$, the time of $(m+1)$-st switching occur after time $t$ can be defined recursively as the first switching after the $m$-th switching after time $t$, that is,
\[
\tau_t^{(m+1)}:=\inf\left\{s>\tau^{(m)}_t:R_s\neq R_{\tau^{(m)}_t}\right\}=\tau_{\tau_t^{(m)}}^{(1)}.
\]
with $\tau^{(0)}_t:=t$.
Let us remind that the stopping times $\tau_t^{(m)}$ here are defined for the regime process $\{R_t:\,t\geq 0\}$ alone, whereas the stopping time $\eta_t$ defined in \eqref{etadef} is defined for the entire stochastic process $\{X_t:t\geq 0\}$.

Hereafter, we reserve $T\in(0,+\infty)$ for the terminal time.
For $i\in\mathcal{M}$, we define the differential operator $\mathcal{L}_i$ by:
\[
\mathcal{L}_{i} f(t,{\bf x}) := \langle \nabla f(t,{\bf x}), b(t,{\bf x},i)\rangle + \frac{1}{2} {\rm tr}\left[\sigma(t,{\bf x},i)^{\top}({\rm Hess}f(t,{\bf x}))\sigma(t,{\bf x},i)\right], %\quad (t,{\bf x})\in[0,T)\times D,
\]
provided that the right-hand side is finite valued, where the gradient and Hessian are taken with respect to the state variable.
Moreover, we reserve $r$ for a bounded continuous real-valued function on $[0,T]\times D \times\mathcal{M}$, and define
\[
\Theta_{t_1,t_2}:=\exp\left[-\int_{t_1}^{t_2}r(s,X_s,R_s) ds\right],\quad \Theta_{t_1,t_2}^i:=\exp\left[-\int_{t_1}^{t_2}r(s,X_s,i) ds\right],\quad \Lambda_{t_1,t_2}^i:=\exp\left[\int_{t_1}^{t_2} q_{ii}(X_s)ds\right],
\]
for $0\leq t_1\leq t_2\leq T$ and $i\in\mathcal{M}$.
Note that the notation $\Theta_{t_1,t_2}^i$ is, strictly speaking, not necessary in some instances (for instance, $\Theta_{t_1,t_2}=\Theta_{t_1,t_2}^i$ for $t\le t_1\le t_2\le T$ under $\mathbb{P}_0^{t,{\bf x},i}$), whereas we employ it for clearer presentation as well as some ease of coding. 
Then, for $(i,\gamma)\in\mathcal{M}\times[0,1]$, $f(\cdot,\cdot,i) \in \mathcal{C}^{1,2}_0([0,T)\times D;\mathbb{R})$ and a generator matrix $Q$, the infinitesimal generator of the paired process $\{(X_s,R_s):s\geq 0\}$ is given by \cite[Theorem 3]{BYZ}:
\begin{equation}\label{gamma}
\lim_{\Delta\to 0}\frac{\mathbb{E}_{\gamma Q}^{t,{\bf x},i}\left[\Theta_{t,t+\Delta}f(t+\Delta,X_{t+\Delta},R_{t+\Delta})\right]-f(t,{\bf x},i)}{\Delta}
=\frac{\partial}{\partial t}f(t,{\bf x},i)+\mathcal{L}_{i}f(t,{\bf x},i)-r(t,{\bf x},i)f(t,{\bf x},i)+\gamma \sum_{j\in\mathcal{M}}q_{ij}({\bf x}) f(t,{\bf x},j),
\end{equation}
where we denote by $\mathcal{C}_0^{1,2}$ the class of functions, whose first derivative in the time variable and second derivatives in the state variables are continuous with compact support.
Since $Q$ is a generator matrix satisfying the $q$-property (Assumption \ref{qproperty}), so is the matrix $\gamma Q(\cdot)$ for any $\gamma \in[0,1]$.
If $\gamma>0$, the infinitesimal generator \eqref{gamma} is said to be weakly coupled, in the sense that the value at $(t,{\bf x},i)$ depends on $f(t,{\bf x},j)$ for $j\ne i$.
If $\gamma=0$, however, such inter-dependency gets untangled, that is, the infinitesimal generator \eqref{gamma} reduces to the usual one $((\partial/\partial t)+ \mathcal{L}_i-r(t,{\bf x},i)) f(t,{\bf x},i)$ in the absence of involvement of $f(\cdot,\cdot,j)$ for $j\ne i$.
%hich is independent of $f(t,{\bf x},j)$, and
%$\mathbb{E}_0^{t,{\bf x},i}$ denotes the expectation under $\mathbb{P}_0^{t,{\bf x},i}$,
%with respect to the generator matrix $Q(\cdot)\equiv 0$.
If $\gamma=0$, then the corresponding generator matrix is the zero matrix in $\mathbb{R}^{p\times p}$.
Hence, the probability of a switching occur within an infinitesimal-sized interval $(t,t+\Delta]$ is as negligibly small as $o(\Delta)$ due to \eqref{defPR}.

\section{Initial value problems}\label{section initial value problems}

We are concerned with weak approximation of the following function written in terms of mathematical expectation:
\begin{equation}\label{defv}
v(t,{\bf x},i):= \mathbb{E}_Q^{t,{\bf x},i} \left[ \Theta_{t,T}  g(X_T, R_T) - 
\int_t^T \Theta_{t,s}\phi(s,X_s,R_s)ds \right], \quad (t,{\bf x},i) \in  [0,T]\times D \times \mathcal{M},
\end{equation}
with $g:D\times \mathcal{M}\to\mathbb{R}$ and $\phi:[0,T]\times D\times \mathcal{M}\to \mathbb{R}$.
For instance, Monte Carlo approximation of the probabilistic representation \eqref{defv} is not trivial in the presence of switching ($Q\ne 0$), clearly because one needs to keep track of the additional regime process $\{R_s:\,s\in [t,T]\}$ as well as modulate the diffusion process $\{X_s:\,s\in [t,T]\}$ according to the switching. 
It is known \cite{BYZ, doi:10.1142/p473, YZ1} that the probabilistic representation \eqref{defv} is a unique solution to the following initial value problem:
\begin{equation}\label{pdev}
\begin{dcases}
\frac{\partial}{\partial t}v(t,{\bf x},i)+\mathcal{L}_{i}v( t,{\bf x},i) =  (r(t,{\bf x},i)-q_{ii}({\bf x}))v(t,{\bf x},i)+\phi(t,{\bf x},i)-\sum_{j\in\mathcal{M}\setminus \{i\}}q_{ij}({\bf x}) v(t,{\bf x},j),
&  (t,{\bf x},i) \in [0,T)\times D  \times \mathcal{M},\\
 v(T,{\bf x},i)=g({\bf x},i), & ({\bf x},i) \in D \times \mathcal{M},
\end{dcases}
\end{equation}
under the following regularity conditions on the input data.

\begin{assumption}\label{Ai}
There exists an $\alpha\in(0,1)$, such that for every $i\in\mathcal{M}$,

\noindent
{\bf (a)}  $r(\cdot,\cdot,i)$ is bounded on $[0,T]\times D$, and is $\alpha$-H\"older continuous on every compact subset of $[0,T]\times D$,

\noindent
{\bf (b)}   $\phi(\cdot,{\bf x},i)$ is continuous on $[0,T]$ for all ${\bf x}\in D$,

\noindent 
{\bf (c)} $\phi(t,\cdot,i)$ is $\alpha$-H\"older continuous and at most of polynomial growth on $D$ for all $t\in [0,T]$,

\noindent
{\bf (d)}  $g(\cdot,i)$ is continuous and at most of polynomial growth on $D$.
\end{assumption}

The translation of the probabilistic representation \eqref{defv} into the initial value problem \eqref{pdev} does not allow us to avoid the additional complexity caused by the switching, because here $p(=|\mathcal{M}|)$ partial differential equations in \eqref{pdev} are not mutually independent but weakly coupled due to the last term $\sum_{j\in\mathcal{M}\setminus \{i\}}q_{ij}({\bf x}) v(t,{\bf x},j)$.
%in fact, hidden in the differential operator \eqref{operator Ai}.
The primary contribution of this work is to construct an iterative weak approximation framework (Sections \ref{section iterative weak approximation} and \ref{section monotonic}) and associated hard upper and lower bounding functions (Section \ref{section hard bounding functions IVP}) with suitable convergence towards the target solution \eqref{defv}, without dealing with such regime-switching or inter-dependency of multiple partial differential equations.

\subsection{Iterative weak approximation}\label{section iterative weak approximation}
In order to construct our iterative weak approximation framework, we begin with the following function on $[0,T]\times D\times \mathcal{M}$ in terms of mathematical expectation:
\begin{equation}\label{def of w0}
w_0(t,{\bf x},i):= \mathbb{E}_0^{t,{\bf x},i} \left[ \Theta_{t,T}^i  g(X_T,i) - 
\int_t^T \Theta_{t,s}^i\phi(s,X_s,i)ds \right], \quad (t,{\bf x},i) \in  [0,T]\times D \times \mathcal{M},
\end{equation}
where the expectation is taken under the probability measure $\mathbb{P}_0^{t,{\bf x},i}$ with zero generator matrix ($Q\equiv 0$), meaning that the regime process cannot switch its regime but remains at the initial regime throughout (that is, $R_s\equiv i$ for all $s\in [t,T]$) under this probability measure.
Note that we specify the superscript $i$ of the notation $\Theta_{t,\cdot}^i$ in \eqref{def of w0} for completeness and slight ease of coding, although the superscript can be suppressed without confusion.  
Hence, under Assumption \ref{Ai}, the function \eqref{def of w0} uniquely solves the initial value problem:
\begin{equation}\label{IVPw0pde}
\begin{dcases}
\frac{\partial}{\partial t}w_0(t,{\bf x},i)+\mathcal{L}_{i}w_0( t,{\bf x},i)  =  r(t,{\bf x},i)w_0(t,{\bf x},i)+\phi(t,{\bf x},i),
&  (t,{\bf x},i) \in [0,T)\times D  \times \mathcal{M},\\
 w_0(T,{\bf x},i)=g({\bf x},i),
&  ({\bf x},i) \in  D  \times \mathcal{M}.
\end{dcases}
\end{equation}
Next, with $w_0$ as an initial guess, we construct the sequence $\{w_m\}_{m\in\mathbb{N}}$ of continuous functions by recursion:
\begin{equation}\label{IVPwm}
w_m(t,{\bf x},i):=\mathbb{E}_0^{t,{\bf x},i}\left[\Theta_{t,T}^i\Lambda^i_{t,T}g(X_T,i) -\int_t^T \Theta_{t,s}^i\Lambda^i_{t,s} \left(\phi(s,X_s,i) -\sum_{j\in\mathcal{M}\setminus\{i\}}q_{ij}(X_s)w_{m-1}(s,X_s,j) \right)ds\right],
\end{equation}
which uniquely solves the initial value problem:
\begin{equation}\label{IVPwmpde}
\begin{dcases}
 \frac{\partial}{\partial t} w_m(t,{\bf x},i) +\mathcal{L}_i w_m(t,{\bf x},i) = (r(t,{\bf x},i)-q_{ii}({\bf x}))w_m(t,{\bf x},i)+\phi(t,{\bf x},i)-\sum_{j\in\mathcal{M}\setminus\{i\}}q_{ij}({\bf x})w_{m-1}(t,{\bf x},j),
& (t,{\bf x},i) \in [0,T)\times D \times \mathcal{M},\\
w_m(T,{\bf x},i) =g({\bf x},i),
& ({\bf x},i) \in D \times \mathcal{M},
\end{dcases}
\end{equation}
provided that the term $\phi(t,{\bf x},i)-\sum_{j\in\mathcal{M}\setminus\{i\}}q_{ij}({\bf x})w_{m-1}(t,{\bf x},j)$ is smooth enough (Assumptions \ref{qproperty} and \ref{Ai}) to be treated as a heat source on the whole. 

%As a numerical method, the proposed iteration significantly reduces the complexity caused by the switching and the inter-dependency among partial differential equations, because the expectation and are taken under the probability measure.

The iteration \eqref{def of w0}-\eqref{IVPwmpde} significantly reduces the complexity caused by the switching \eqref{defv} and the inter-dependency among partial differential equations \eqref{pdev}, because the expectations \eqref{def of w0} and \eqref{IVPwm} are taken under the probability measure $\mathbb{P}_0^{t,{\bf x},i}$, that is, no switching is possible due to the zero generator matrix ($Q(\cdot)\equiv 0$), whereas the initial value problems \eqref{IVPw0pde} and \eqref{IVPwmpde} are not weakly coupled but only a collection of $p(=|\mathcal{M}|)$ independent initial value problems, based on the given the previous iterate $w_{m-1}$ in \eqref{IVPwmpde}.
Hence, each step only entails the computation of the standard (that is, non-switching) diffusion process or of a linear initial value problem, for which a variety of numerical methods are available, such as discretization methods for (non-switching) standard stochastic differential equations, finite difference and element methods, to name just a few.  

A natural and crucial theoretical question on such iterative schemes is whether or not the iteration under consideration converges to the desired solution.
We are now in a position to give the main result of this section, that is, we prove positive.
%that our iterative built a theoretical base for a convergent iterative weak approximation method.
%We prove rigorously that our iteration indeed does.
%Unlike the usual proofs via partial differential equations, we employ a probabilistic approach.

%\textcolor{red}{
%Precisely, we prove that any two consecutive functions of \eqref{wmdef} (that is, $w_{m-1}$ and $w_m$ for any $m\in\mathbb{N}$) can be connected in the form of the probabilistic representation \eqref{IVPwm} or equivalently the initial value problem \eqref{IVPwmpde}, as well as that the sequence $\{w_m\}_{m\in\mathbb{N}_0}$ converges to the target function \eqref{defv} pointwise, or even locally uniformly under an additional condition.}
%We recall the notation $\Theta_{t_1,t_2}^i$ and $\Lambda_{t_1,t_2}^i$ defined in \eqref{def of thetai}.

\begin{theorem}\label{IVPwtheorem}
Suppose that the functions $w_m(\cdot,\cdot,i)$ defined in \eqref{def of w0} and \eqref{IVPwm} are in $\mathcal{C}_0^{1,2}([0,T)\times D;\,\mathbb{R})\cap \mathcal{C}([0,T]\times D;\,\mathbb{R})$ for all $m\in\mathbb{N}_0$ and $i\in\mathcal{M}$.

\noindent
{\bf (a)}
It holds that for $m\in\mathbb{N}$ and $(t,{\bf x},i) \in [0,T]\times D \times\mathcal{M}$,
\begin{equation}\label{wmdef}
w_m(t,{\bf x},i)=\mathbb{E}_Q^{t,{\bf x},i}\left[\Theta_{t,\tau_t^{(m)}\land T}w_0(\tau_t^{(m)}\land T,X_{\tau_t^{(m)}\land T},R_{\tau_t^{(m)}\land T}) - \int_t^{\tau_t^{(m)} \land T} \Theta_{t,s}\phi(s,X_s,R_s)ds\right].
\end{equation}

\noindent
{\bf (b)}
It holds that for $(t,{\bf x},i) \in [0,T]\times D \times\mathcal{M}$, $|w_m(t,{\bf x}, i)-v(t,{\bf x},i)|=\mathcal{O}((c^m(T-t)^m/m!)^p)$ as $m\to +\infty$ for any $p\in (0,1)$, with $c:=\sup_{{\bf x}\in D,j\ne k}q_{jk}({\bf x})$.
Moreover, if the functions $\{w_0(\cdot,\cdot,i)\}_{i\in\mathcal{M}}$ are either all non-negative or all non-positive on $[0,T]\times D$, then we have $w_m(\cdot, \cdot,i)\to v(\cdot,\cdot,i)$ locally uniformly on $[0,T]\times D$.
\end{theorem}

%To maintain the flow of the paper, we defer the proof to Appendix \ref{section proofs}.
In summary, we prove (Theorem \ref{IVPwtheorem} {\bf (a)}) that each iterate $w_m$, defined recursively in \eqref{IVPwm} {\it without switching}, represents the expression \eqref{wmdef} {\it with switching}, in which the underlying dynamics is terminated when a given number of switching occur before the terminal time, that is, $\tau_t^{(m)}\land T$.
That is, as is observable from the representation \eqref{wmdef}, as the restriction gets relaxed ($m\to +\infty$), the sequence $\{w_m\}_{m\in\mathbb{N}}$ tends to the target function \eqref{defv} as desired, due to $\tau_t^{(m)}\to +\infty$ by Assumption \ref{qproperty} {\bf (c)}. 
Note that the result \eqref{wmdef} does not recover the initial guess \eqref{def of w0} but still remains consistent with $m=0$, in the sense that \eqref{wmdef} yields the trivial identity $w_0(t,{\bf x},i)=w_0(t,{\bf x},i)$, as we have set $\tau^{(0)}_t=t$ for all $t\in [0,T]$.
Hence, with theoretical guarantee of convergence and its rate in Theorem \ref{IVPwtheorem} {\bf (b)}, we repeat this iteration until we somehow decide to terminate it.

\begin{proof}[Proof of Theorem \ref{IVPwtheorem}]\label{proofA1}
{\bf (a)}
In light of Assumptions \ref{qproperty} and \ref{Ai}, the function $q_{ii}(\cdot)$ is bounded and as smooth as $r(t,\cdot,i)$.
Moreover, the function $\sum_{j\in\mathcal{M}\setminus\{i\}}q_{ij}(\cdot)w_{m-1}(\cdot,\cdot,j)$ is at least as smooth as $\phi(\cdot,\cdot,i)$, for all $i\in\mathcal{M}$, that is, the probabilistic representation \eqref{IVPwm} is a unique solution to the initial value problem \eqref{IVPwmpde}.
It thus suffices to show that the desired representation \eqref{wmdef} solves the initial value problem \eqref{IVPwmpde} for all $m\in\mathbb{N}$.
To this end, for $t\in [0,T]$ and $m\in\mathbb{N}$, define 
\[
F_{t,m} := \Theta_{t,\tau_t^{(1)}\land T}w_{m-1}(\tau_t^{(1)}\land T,X_{\tau_t^{(1)}\land T},R_{\tau_t^{(1)}\land T}) 
- \int_t^{\tau_t^{(1)} \land T} \Theta_{t,s}\phi(s,X_s,R_s)ds.
\]
We first show that for $m\in\mathbb{N}$, the representation \eqref{wmdef} can be rewritten as $w_m(t,{\bf x},i)=\mathbb{E}_Q^{t,{\bf x},i}[F_{t,m}]$.
The case $m=1$ is trivial as it simply reduces to the representation \eqref{wmdef}.
%\[
%w_1(t,{\bf x},i)=\mathbb{E}_Q^{t,{\bf x},i}\left[\Theta_{t,\tau_t^{(1)}\land T}w_{0}(\tau_t^{(1)}\land T,X_{\tau_t^{(1)}\land T},R_{\tau_t^{(1)}\land T}) - \int_t^{\tau_t^{(1)} \land T} \Theta_{t,s}\phi(s,X_s,R_s)ds\right].
%\]
Let $m\in\{2,\cdots\}$.
It holds that for $(t,{\bf x},i) \in  [0,T]\times D \times \mathcal{M}$, 
\begin{align}
w_m(t,{\bf x},i){}&=\mathbb{E}_Q^{t,{\bf x},i}\left[\Theta_{t,\tau_t^{(m)}\land T}w_{0}(\tau_t^{(m)}\land T,X_{\tau_t^{(m)}\land T},R_{\tau_t^{(m)}\land T}) 
- \int_t^{\tau_t^{(m)} \land T} \Theta_{t,s}\phi(s,X_s,R_s)ds\right] \nonumber\\
%&=\mathbb{E}_Q^{t,{\bf x},i}\left[\mathbbm{1}(\tau_t^{(1)}>T) \left(\Theta_{t,\tau_t^{(m)}\land T}w_{0}(\tau_t^{(m)}\land T,X_{\tau_t^{(m)}\land T},R_{\tau_t^{(m)}\land T}) 
%- \int_t^{\tau_t^{(m)} \land T} \Theta_{t,s}\phi(s,X_s,R_s)ds\right)\right] \nonumber\\
%&\qquad +\mathbb{E}_Q^{t,{\bf x},i}\left[\mathbbm{1}(\tau_t^{(1)}\leq T) \left(\Theta_{t,\tau_t^{(m)}\land T}w_{0}(\tau_t^{(m)}\land T,X_{\tau_t^{(m)}\land T},R_{\tau_t^{(m)}\land T}) 
%- \int_t^{\tau_t^{(m)} \land T} \Theta_{t,s}\phi(s,X_s,R_s)ds\right)\right] \nonumber\\
&=\mathbb{E}_Q^{t,{\bf x},i}\Bigg[\mathbbm{1}(\tau_t^{(1)}>T) \left(\Theta_{t,T}w_0(T, X_T,R_T)- \int_t^{ T} \Theta_{t,s}\phi(s,X_s,R_s)ds\right) -\mathbbm{1}(\tau_t^{(1)}\leq T) \int_t^{\tau_t^{(1)}} \Theta_{t,s}\phi(s,X_s,R_s)ds \nonumber\\
&\qquad\qquad  +\mathbbm{1}(\tau_t^{(1)}\leq T)\Theta_{t,\tau_t^{(1)}} \left(\Theta_{\tau_t^{(1)},\tau_t^{(m)}\land T}w_{0}(\tau_t^{(m)}\land T,X_{\tau_t^{(m)}\land T},R_{\tau_t^{(m)}\land T})   -\int_{\tau_t^{(1)}}^{\tau_t^{(m)}\land T}\Theta_{\tau_t^{(1)},s}\phi(s,X_s,R_s)ds\right)\Bigg], \label{lemmapA1}
\end{align}
where the last equality holds because $\tau_t^{(m)}\land T=T$ on the event $\{\tau_t^{(1)}>T\}$.
Then, it holds that
\begin{align}
&\mathbb{E}_Q^{t,{\bf x},i}\left[\mathbbm{1}(\tau_t^{(1)}\leq T)\Theta_{t,\tau_t^{(1)}} \left(\Theta_{\tau_t^{(1)},\tau_t^{(m)}\land T}w_{0}(\tau_t^{(m)}\land T,X_{\tau_t^{(m)}\land T},R_{\tau_t^{(m)}\land T})   -\int_{\tau_t^{(1)}}^{\tau_t^{(m)}\land T}\Theta_{\tau_t^{(1)},s}\phi(s,X_s,R_s)ds\right)\right] \nonumber \\
&\, =\mathbb{E}_Q^{t,{\bf x},i}\left[\mathbbm{1}(\tau_t^{(1)}\leq T)\Theta_{t,\tau_t^{(1)}} \left(\Theta_{\tau_t^{(1)},\tau_{\tau_t^{(1)}}^{(m-1)}\land T}w_{0}(\tau_{\tau_t^{(1)}}^{(m-1)}\land T,X_{\tau_{\tau_t^{(1)}}^{(m-1)}\land T},R_{\tau_{\tau_t^{(1)}}^{(m-1)}\land T})   -\int_{\tau_t^{(1)}}^{\tau_{\tau_t^{(1)}}^{(m-1)}\land T}\Theta_{\tau_t^{(1)},s}\phi(s,X_s,R_s)ds\right)\right] \nonumber \\
&\, =\mathbb{E}_Q^{t,{\bf x},i}\left[\mathbbm{1}(\tau_t^{(1)}\leq T)\Theta_{t,\tau_t^{(1)}} \mathbb{E}_Q^{t,{\bf x},i}\left[\Theta_{\tau_t^{(1)},\tau_{\tau_t^{(1)}}^{(m-1)}\land T}w_{0}(\tau_{\tau_t^{(1)}}^{(m-1)}\land T,X_{\tau_{\tau_t^{(1)}}^{(m-1)}\land T},R_{\tau_{\tau_t^{(1)}}^{(m-1)}\land T})   -\int_{\tau_t^{(1)}}^{\tau_{\tau_t^{(1)}}^{(m-1)}\land T}\Theta_{\tau_t^{(1)},s}\phi(s,X_s,R_s)ds\Bigg|\, \mathcal{F}_{\tau_t^{(1)}}\right]\right]\nonumber \\
&\, =\mathbb{E}_Q^{t,{\bf x},i}\left[\mathbbm{1}(\tau_t^{(1)}\leq T)\Theta_{t,\tau_t^{(1)}} \mathbb{E}_Q^{\tau_t^{(1)},X_{\tau_t^{(1)}},R_{\tau_t^{(1)}}}\left[\Theta_{\tau_t^{(1)},\tau_{\tau_t^{(1)}}^{(m-1)}\land T}w_{0}(\tau_{\tau_t^{(1)}}^{(m-1)}\land T,X_{\tau_{\tau_t^{(1)}}^{(m-1)}\land T},R_{\tau_{\tau_t^{(1)}}^{(m-1)}\land T})   -\int_{\tau_t^{(1)}}^{\tau_{\tau_t^{(1)}}^{(m-1)}\land T}\Theta_{\tau_t^{(1)},s}\phi(s,X_s,R_s)ds\right]\right]\nonumber \\
&\, =\mathbb{E}_Q^{t,{\bf x},i}\left[\mathbbm{1}(\tau_t^{(1)}\leq T)\Theta_{t,\tau_t^{(1)}} w_{m-1}(\tau_t^{(1)},X_{\tau_t^{(1)}},R_{\tau_t^{(1)}})\right], \label{lemmapA2}
\end{align}
where the first equality holds by $\tau_t^{(m)}=\tau_{\tau_t^{(1)}}^{(m-1)}$ for all $m\in\mathbb{N}$, the second by the tower property, the third by the strong Markov property and the last one by the representation \eqref{wmdef}.
Here, we denote by $\mathcal{F}_{\tau_t^{(1)}}$ the stopping time $\sigma$-field at $\tau_t^{(1)}$, that is, $\mathcal{F}_{\tau_t^{(1)}}:=\{B\in \mathcal{F}:\,B\cap \{\tau_t^{(1)}\le s\}\in \mathcal{F}_s \text{ for all }s \in (t,T]\}.$
On the whole, combining \eqref{lemmapA1} and \eqref{lemmapA2} yields the representation $w_m(t,{\bf x},i)=\mathbb{E}_Q^{t,{\bf x},i}[F_{t,m}]$ for all $m\in\mathbb{N}$, as desired.
%\begin{align*}
%w_m(t,{\bf x},i)&=\mathbb{E}_Q^{t,{\bf x},i}\Bigg[\mathbbm{1}(\tau_t^{(1)}>T) \Theta_{t,T}w_{m-1}(T, X_T,R_T)+\mathbbm{1}(\tau_t^{(1)}\leq T)\Theta_{t,\tau_t^{(1)}} w_{m-1}(\tau_t^{(1)},X_{\tau_t^{(1)}},R_{\tau_t^{(1)}})  \\
%&\qquad \qquad -  \mathbbm{1}(\tau_t^{(1)}> T)\int_{t}^{T}\Theta_{t,s}\phi(s,X_s,R_s)ds - \mathbbm{1}(\tau_t^{(1)}\leq T) \int_t^{\tau_t^{(1)}} \Theta_{t,s}\phi(s,X_s,R_s)ds \Bigg] \\
%&=\mathbb{E}_Q^{t,{\bf x},i}\left[\Theta_{t,\tau_t^{(1)}\land T}w_{m-1}(\tau_t^{(1)}\land T,X_{\tau_t^{(1)}\land T},R_{\tau_t^{(1)}\land T}) 
%- \int_t^{\tau_t^{(1)} \land T} \Theta_{t,s}\phi(s,X_s,R_s)ds\right],
%\end{align*}
%which concludes.

Next, we split the expectation $\mathbb{E}_Q^{t,{\bf x},i}[F_{t,m}]$ in terms of the destination of the regime process after the first switching.
To this end, fix $\delta\in (0,T-t]$ and denote by $J_{t,t+\delta}^i :=\{\tau_t^{(1)}>t+\delta\}$ the event of no switching within the interval $(t,t+\delta]$, and denote by $J_{t,t+\delta}^j:=\{\tau_t^{(1)}\leq t+\delta\}\cap \{R_{\tau_t^{(1)}}=j\}$ the event that the first switching occur within the interval $(t,t+\delta]$ and the destination is regime $j$.
With $(t,i)\in [0,T]\times \mathcal{M}$ and $\delta\in (0,T-t]$ fixed, the events $\{J_{t,t+\delta}^j\}_{j\in\mathcal{M}}$ are clearly jointly exhaustive, that is, $\mathbbm{1}(J^i_{t,t+\delta})+\sum_{j\in \mathcal{M}\setminus \{i\}}\mathbbm{1}(J^j_{t,t+\delta})=1$, and thus \color{black}
\begin{equation}\label{IVTp1}
w_m(t,{\bf x},i)=\mathbb{E}_Q^{t,{\bf x},i}\left[F_{t,m}\right] 
= \mathbb{E}_Q^{t,{\bf x},i}\left[\mathbbm{1}(J_{t,t+\delta}^i) F_{t,m}\right] +\sum_{j\in\mathcal{M}\setminus\{i\}}\mathbb{E}_Q^{t,{\bf x},i}\left[\mathbbm{1}(J_{t,t+\delta}^j) F_{t,m}\right],\quad (t,{\bf x},i) \in  [0,T]\times D \times \mathcal{M}.
\end{equation}
%We investigate the rightmost representation of the result \eqref{IVTp1}.
As for the first expectation $\mathbb{E}_Q^{t,{\bf x},i}[\mathbbm{1}(J_{t,t+\delta}^i) F_{t,m}]$, the event $J_{t,t+\delta}^i$ ensures no switching within the interval $(t,t+\delta]$.
% that is, the first switching after time $t$, if at all, can occur within the remaining interval $(t+\delta,T]$.
Hence, we have $\tau_t^{(1)} = \tau_{t+\delta}^{(1)}$ on the event $J_{t,t+\delta}^i$, and moreover, 
\begin{align}
\mathbb{E}_Q^{t,{\bf x},i}\left[\mathbbm{1}(J_{t,t+\delta}^i) F_{t,m}\right] 
&=\mathbb{E}_Q^{t,{\bf x},i}\left[\mathbbm{1}(J_{t,t+\delta}^i) \left(\Theta_{t,\tau_t^{(1)}\land T}w_{m-1}(\tau_t^{(1)}\land T,X_{\tau_t^{(1)}\land T},R_{\tau_t^{(1)}\land T}) 
- \int_t^{\tau_t^{(1)} \land T} \Theta_{t,s}\phi(s,X_s,R_s)ds\right)\right] \nonumber \\
%&=\mathbb{E}_Q^{t,{\bf x},i}\Bigg[\mathbbm{1}(J_{t,t+\delta}^i) \Theta_{t,t+\delta} \left(\Theta_{t+\delta,\tau_{t+\delta}^{(1)}\land T}w_{m-1}(\tau_{t+\delta}^{(1)}\land T,X_{\tau_{t+\delta}^{(1)}\land T},R_{\tau_{t+\delta}^{(1)}\land T}) - \int_{t+\delta}^{\tau_{t+\delta}^{(1)} \land T} \Theta_{t+\delta,s}\phi(s,X_s,R_s)ds\right)\nonumber \\
%&\qquad\qquad\qquad  - \mathbbm{1}(J_{t,t+\delta}^i) \int_t^{t+\delta} \Theta_{t,s}\phi(s,X_s,R_s)ds\Bigg] \nonumber\\
&=\mathbb{E}_Q^{t,{\bf x},i}\left[\mathbbm{1}(J_{t,t+\delta}^i) \Theta_{t,t+\delta} F_{t+\delta,m}-\mathbbm{1}(J_{t,t+\delta}^i) \int_t^{t+\delta} \Theta_{t,s}\phi(s,X_s,R_s)ds\right] \nonumber\\
&=\mathbb{E}_Q^{t,{\bf x},i}\left[\mathbbm{1}(J_{t,t+\delta}^i) \Theta_{t,t+\delta} \mathbb{E}_Q^{t+\delta,X_{t+\delta},R_{t+\delta}}\left[ F_{t+\delta,m}\right]-\mathbbm{1}(J_{t,t+\delta}^i) \int_t^{t+\delta} \Theta_{t,s}\phi(s,X_s,R_s)ds\right] \nonumber\\
&=\mathbb{E}_Q^{t,{\bf x},i}\left[\mathbbm{1}(J_{t,t+\delta}^i) \Theta_{t,t+\delta} w_m(t+\delta,X_{t+\delta},R_{t+\delta})-\mathbbm{1}(J_{t,t+\delta}^i) \int_t^{t+\delta} \Theta_{t,s}\phi(s,X_s,R_s)ds\right], \label{IVTp2}
\end{align}
where we have applied the tower and Markov properties for the third equality and the result \eqref{IVTp1} for the last equality.
As for the second expectation $\mathbb{E}_Q^{t,{\bf x},i}[\mathbbm{1}(J_{t,t+\delta}^j) F_{t,m}]$, the event $J_{t,t+\delta}^j$ ensures switching at least once within the interval $(t,t+\delta]$, that is, $\tau_t^{(1)}\leq t+\delta\leq T$.
Hence, we have 
\begin{align}
\mathbb{E}_Q^{t,{\bf x},i}\left[\mathbbm{1}(J_{t,t+\delta}^j) F_{t,m}\right] 
&=\mathbb{E}_Q^{t,{\bf x},i}\left[\mathbbm{1}(J_{t,t+\delta}^j) \left(\Theta_{t,\tau_t^{(1)}\land T}w_{m-1}(\tau_t^{(1)}\land T,X_{\tau_t^{(1)}\land T},R_{\tau_t^{(1)}\land T}) 
- \int_t^{\tau_t^{(1)} \land T} \Theta_{t,s}\phi(s,X_s,R_s)ds\right)\right] \nonumber  \\
&=\mathbb{E}_Q^{t,{\bf x},i}\left[\mathbbm{1}(J_{t,t+\delta}^j) \Theta_{t,\tau_t^{(1)}}w_{m-1}(\tau_t^{(1)},X_{\tau_t^{(1)}},R_{\tau_t^{(1)}})-\mathbbm{1}(J_{t,t+\delta}^j)\int_t^{\tau_t^{(1)}} \Theta_{t,s}\phi(s,X_s,R_s)ds\right]. \label{IVTp3}
\end{align}
Finally, combining (\ref{IVTp1}), (\ref{IVTp2}) and (\ref{IVTp3}) yields
\begin{multline}\label{IVTp31}
w_m(t,{\bf x},i)\\
=\mathbb{E}_Q^{t,{\bf x},i}\left[\mathbbm{1}(J_{t,t+\delta}^i) \Theta_{t,t+\delta} w_m(t+\delta,X_{t+\delta},i)+\sum_{j\in\mathcal{M}\setminus\{i\}}\mathbbm{1}(J_{t,t+\delta}^j)\Theta_{t,\tau_t^{(1)}}w_{m-1}(\tau_t^{(1)},X_{\tau_t^{(1)}},j)-\int_t^{\tau_t^{(1)}\land (t+\delta)} \Theta_{t,s}\phi(s,X_s,R_s)ds\right],
\end{multline}
for all $\delta\in[0,T-t]$, due to $\mathbbm{1}(J^i_{t,t+\delta}) = \mathbbm{1}(\tau_t^{(1)}>t+\delta)$ and $\sum_{j\in\mathcal{M}\setminus\{i\}}\mathbbm{1}(J^j_{t,t+\delta}) = \mathbbm{1}(\tau_t^{(1)}\leq t+\delta)$.
With the function $f_i:[0,T]\times D\times\mathcal{M}\to  \mathbb{R}$ defined by
\[
f_i(t,{\bf x},j):=
\begin{dcases}
w_m(t,{\bf x},i), &\text{if } j=i,\\
w_{m-1}(t,{\bf x},j), &\text{if } j\in \mathcal{M}\setminus \{i\},
\end{dcases}
\]
%which satisfies the following properties:
%\[
%\begin{dcases}
%f_i(t+\delta,X_{t+\delta},R_{t+\delta}) = \mathbbm{1}(R_{t+\delta}=i) w_m(t+\delta,X_{t+\delta},i) +\mathbbm{1}(R_{t+\delta} \neq i) w_{m-1}(t+\delta,X_{t+\delta},R_{t+\delta}),\\
%\mathbbm{1}(J_{t,t+\delta}^i) f_i(t+\delta,X_{t+\delta},R_{t+\delta}) = \mathbbm{1}(J_{t,t+\delta}^i) w_m(t+\delta,X_{t+\delta},i), \\
%\mathbbm{1}(J_{t,t+\delta}^j) f_i(\tau_t^{(1)},X_{\tau_t^{(1)}},R_{\tau_t^{(1)}}) = \mathbbm{1}(J_{t,t+\delta}^j) w_{m-1}(\tau_t^{(1)},X_{\tau_t^{(1)}},j), \quad \text{for } j\neq i. 
%\end{dcases}
%\]
the identity (\ref{IVTp31}) can be rewritten as %, we obtain that for all $\delta\in[0,T-t]$,
\begin{align*}
&f_i(t,{\bf x},i)\\
&\, =\mathbb{E}_Q^{t,{\bf x},i}\left[\mathbbm{1}(J_{t,t+\delta}^i) \Theta_{t,t+\delta} f_i(t+\delta,X_{t+\delta},R_{t+\delta})+\sum_{j\in\mathcal{M}\setminus\{i\}}\mathbbm{1}(J_{t,t+\delta}^j)\Theta_{t,\tau_t^{(1)}}f_i(\tau_t^{(1)},X_{\tau_t^{(1)}},R_{\tau_t^{(1)}})-\int_t^{\tau_t^{(1)}\land (t+\delta)} \Theta_{t,s}\phi(s,X_s,R_s)ds\right] \nonumber\\
&\,=\mathbb{E}_Q^{t,{\bf x},i}\left[ \Theta_{t,\tau_t^{(1)}\land (t+\delta)} f_i(\tau_t^{(1)}\land (t+\delta),X_{\tau_t^{(1)}\land (t+\delta)},R_{\tau_t^{(1)}\land (t+\delta)})-\int_t^{\tau_t^{(1)}\land (t+\delta)} \Theta_{t,s}\phi(s,X_s,R_s)ds\right].% \label{IVTp32}
\end{align*}
Note that $f_i(\cdot,\cdot,j)\in\mathcal{C}_0^{1,2}$ for all $j\in\mathcal{M}$, since $w_m(\cdot,\cdot,i)\in\mathcal{C}_0^{1,2}$ and $w_{m-1}(\cdot,\cdot,j)\in\mathcal{C}_0^{1,2}$ for all $j\in\mathcal{M}\setminus\{i\}$. 
By rearranging the last expectation, we obtain
\begin{align*}
0&=\mathbb{E}_Q^{t,{\bf x},i}\left[ \Theta_{t,\tau_t^{(1)}\land (t+\delta)} f_i(\tau_t^{(1)}\land (t+\delta),X_{\tau_t^{(1)}\land (t+\delta)},R_{\tau_t^{(1)}\land (t+\delta)}) \right] - f_i(t,{\bf x},i)-\mathbb{E}_Q^{t,{\bf x},i}\left[\int_t^{\tau_t^{(1)}\land (t+\delta)} \Theta_{t,s}\phi(s,X_s,R_s)ds\right] \nonumber\\
&= \mathbb{E}_Q^{t,{\bf x},i}\left[ \int_t^{\tau_t^{(1)}\land (t+\delta)} \Theta_{t,s} \left[ \left(\frac{\partial}{\partial s} +\mathcal{L}_{R_s} -r(s,X_s,R_s) \right)f_i(s,X_s,R_s) 
+\sum_{j\in\mathcal{M}}q_{R_s j} (X_s) f_i(s,X_s,j) -\phi(s,X_s,R_s)\right] ds\right] \nonumber\\
&= \mathbb{E}_Q^{t,{\bf x},i}\left[ \int_t^{t+\delta } \mathbbm{1}(J_{t,s}^i) \Theta_{t,s} \left[ \left(\frac{\partial}{\partial s} +\mathcal{L}_{i} -r(s,X_s,i) \right)f_i(s,X_s,i) 
+\sum_{j\in\mathcal{M}}q_{i j} (X_s) f_i(s,X_s,j) - \phi(s,X_s,i)\right] ds\right] \nonumber\\
&= \mathbb{E}_Q^{t,{\bf x},i}\left[ \int_t^{t+\delta } \mathbbm{1}(\tau_t^{(1)}>s) \Theta_{t,s} \left[ \left(\frac{\partial}{\partial s} +\mathcal{L}_{i} -r(s,X_s,i) +q_{ii}(X_s)\right)w_m(s,X_s,i) +\sum_{j\in\mathcal{M}\setminus\{i\}}q_{i j} (X_s) w_{m-1}(s,X_s,j) -\phi(s,X_s,i) \right] ds\right],\end{align*}
which holds true for all $\delta\in[0,T-t]$, where the second equality holds by the Dynkin formula with the infinitesimal generator \eqref{gamma}.
Hence, by dividing throughout by $\delta$ and taking the limit $\delta\to  0$, we obtain the desired partial differential equation in the initial value problem \eqref{IVPwmpde}.
Moreover, by setting $t=T$ in \eqref{wmdef}, we obtain
\[
w_m(T,{\bf x},i)= \mathbb{E}_Q^{T,{\bf x},i}\left[\Theta_{T,\tau_T^{(m)}\land T}w_0(\tau_T^{(m)}\land T,X_{\tau_T^{(m)}\land T},R_{\tau_T^{(m)}\land T}) - \int_T^{\tau_T^{(m)} \land T} \Theta_{t,s}\phi(s,X_s,R_s)ds\right] = w_0(T,{\bf x},i) = g({\bf x},i),
\]
due to $\tau_T^{(m)}\land T=T$.
Hence, the representation \eqref{wmdef} satisfies the initial value problem \eqref{IVPwmpde}. 

\noindent {\bf (b)}
We first rewrite \eqref{defv} and \eqref{wmdef} in terms of the function $w_0$ as follows:
\begin{align*}
 v(t,{\bf x},i) - w_m(t,{\bf x},i)&=\mathbb{E}_Q^{t,{\bf x},i}\left[\Theta_{t,T}g(X_T,R_T)-\Theta_{t,\tau_t^{(m)}\land T}w_{0}(\tau_t^{(m)}\land T,X_{\tau_t^{(m)}\land T},R_{\tau_t^{(m)}\land T})- \int_{\tau_t^{(m)} \land T}^T \Theta_{t,s}\phi(s,X_s,R_s)ds\right]\\
 &=\mathbb{E}_Q^{t,{\bf x},i}\Bigg[\Theta_{t,T}w_0(T,X_T,R_T)-\Theta_{t,\tau_t^{(m)}\land T}w_{0}(\tau_t^{(m)}\land T,X_{\tau_t^{(m)}\land T},R_{\tau_t^{(m)}\land T})\\
 &\qquad \qquad - \int_{\tau_t^{(m)} \land T}^T \Theta_{t,s}\left(\frac{\partial}{\partial s}w_0(s,X_s,R_s)+ \mathcal{L}_{R_s}w_0(s,X_s,R_s) -r(s,X_s,R_s)w_0(s,X_s,R_s) \right)ds\Bigg]\\
 &= \mathbb{E}_Q^{t,{\bf x},i}\left[\int_{\tau_t^{(m)}\land T}^T \Theta_{t,s} \sum_{j\in\mathcal{M}}q_{R_s j}(X_s)w_0(s,X_s,j)ds\right],
%&=\mathbb{E}_Q^{t,{\bf x},i}\left[\int_{t}^T \mathbbm{1}(\tau_t^{(m)}\leq s) \Theta_{t,s} \sum_{j\in\mathcal{M}}q_{R_s j}(X_s)w_0(s,X_s,j)ds\right].
\end{align*}
where we have applied the initial value problem \eqref{IVPw0pde} for the second equality and the Ito formula for the third equality.
Then, it holds that for any H\"older conjugates $p$ and $q$,
\begin{align*}
\left|v(t,{\bf x},i) - w_m(t,{\bf x},i)\right|&= \left|\mathbb{E}_Q^{t,{\bf x},i}\left[\mathbbm{1}(\tau_t^{(m)}\le T)\int_{\tau_t^{(m)}\land T}^T \Theta_{t,s} \sum_{j\in\mathcal{M}}q_{R_s j}(X_s)w_0(s,X_s,j)ds\right]\right|\\
&\le \left(\mathbb{P}_Q^{t,{\bf x},i}(\tau_t^{(m)}\le T)\right)^{1/p}(T-t)\mathbb{E}_Q^{t,{\bf x},i}\left[\sup_{s\in [t,T]}\left|\Theta_{t,s} \sum_{j\in\mathcal{M}}q_{R_s j}(X_s)w_0(s,X_s,j)\right|^q\right]^{1/q},
\end{align*}
where the second term in the last line is finite with the aid of \cite[Proposition 2.3]{YZ1}.
Next, letting $N_{t,T}$ denote the number of switching of the regime process on $[t,T]$ and $c:=\sup_{{\bf x}\in D,\,j\ne k}q_{jk}({\bf x})$, it holds that 
\[
 \mathbb{P}_Q^{t,{\bf x},i}(\tau_t^{(m)}\le T)=\sum_{k=m}^{+\infty}\mathbb{P}_Q^{t,{\bf x},i}(N_{t,T}=k)\le \sum_{k=m}^{+\infty}e^{-c(T-t)}\frac{c^k(T-t)^k}{k!}=e^{-a_m}\frac{c^m(T-t)^m}{m!},
\]
where the inequality holds for all $m\in \{\lceil c(T-t)\rceil,\cdots\}$ and $a_m$ is a positive constant in $(0,c(T-t))$ depending on $m$.
This yields the desired decay rate.

%, since the region of integration $(\tau_t^{(m)}\land T,T]$ vanishes eventually as $m\to +\infty$ by Assumption \ref{qproperty} {\bf (c)}.
%Hence, it holds $\mathbb{P}_Q^{t,{\bf x},i}$-$a.s.$ that $\lim_{m\to +\infty}\mathbbm{1}(\tau_t^{(m)}\le s)=0$ for all $s\in (t,T]$.
%Moreover, since $w_0\in\mathcal{C}_0^{1,2}$, and $r(\cdot,\cdot,\cdot)$ , $q_{\cdot j}(\cdot)$ are bounded, we have
%\[
% \left|v(t,{\bf x},i) - w_m(t,{\bf x},i) \right| \leq  \mathbb{E}_Q^{t,{\bf x},i}\left[(T- t) \sup_{s\in[t,T]} \left|\Theta_{t,s} \sum_{j\in\mathcal{M}}q_{R_s j}(X_s)w_0(s,X_s,j) \right|\right] <+\infty.
%\]
%, the dominated convergence theorem justifies the passage to the limit, which yields the pointwise convergence $|v(t,{\bf x},i) - w_m(t,{\bf x},i)|\to  0$ as $m\to +\infty$.
%Therefore, we have $w_m$ converges to $v$ pointwise.

Finally, rewrite the last term as sum of two expectations to construct two sequences $\{v_{m,1}(\cdot,\cdot,i)\}_{m\in\mathbb{N}}$ and $\{v_{m,2}(\cdot,\cdot,i)\}_{m\in\mathbb{N}}$ of real-valued continuous functions for each $i\in\mathcal{M}$, one with the summation over $j\in \mathcal{M}\setminus \{R_s\}$ and the other over $j\in \{R_s\}$, that is,
\[
 v_{m,1}(t,{\bf x},i):=\mathbb{E}_Q^{t,{\bf x},i}\left[\int_{\tau_t^{(m)}\land T}^T \Theta_{t,s} \sum_{j\in\mathcal{M}\setminus\{R_s\}}q_{R_s j}(X_s)w_0(s,X_s,j)ds\right],\quad v_{m,2}(t,{\bf x},i):=\mathbb{E}_Q^{t,{\bf x},i}\left[\int_{\tau_t^{(m)}\land T}^T \Theta_{t,s} q_{R_s R_s}(X_s)w_0(s,X_s,R_s)ds\right],
\]
for $(t,{\bf x},i)\in [0,T]\times D\times \mathcal{M}$.
If the functions $\{w_0(\cdot,\cdot,i)\}_{i\in\mathcal{M}}$ are either all non-negative or all non-positive, then each of the integrands above (that is, $\Theta_{t,s} \sum_{j\in\mathcal{M}\setminus\{R_s\}}q_{R_s j}(X_s)w_0(s,X_s,j)$ and $\Theta_{t,s} q_{R_s R_s}(X_s)w_0(s,X_s,R_s)$) is either non-negative or non-positive for all $s\in (t,T]$, $\mathbb{P}_Q^{t,{\bf x},i}$-$a.s.$
Since the region of integration $(\tau_t^{(m)}\land T,T]$ is non-expanding in $m$, $\mathbb{P}_Q^{t,{\bf x},i}$-$a.s.$, all the sequences $\{v_{m,1}(\cdot,\cdot,i)\}_{m\in\mathbb{N}}$ and $\{v_{m,2}(\cdot,\cdot,i)\}_{m\in\mathbb{N}}$ are monotonic and thus locally uniformly convergent by the Dini theorem.
So are the sums $\{v_{m,1}(\cdot,\cdot,i)+v_{m,2}(\cdot,\cdot,i)\}_{m\in\mathbb{N}}$ for each $i\in\mathcal{M}$, as desired.
\end{proof}

\subsection{Monotonically convergent iterative weak approximation}\label{section monotonic}

We have constructed iterative weak approximations, which are convergent to the unknown solution.
%(Section \ref{section iterative weak approximation}) and associated hard upper and lower bounds (Section \ref{section hard bounding functions IVP}), which are all convergent to the target solution.
In this section, we build an alternative framework on the basis of exactly the same recursion \eqref{IVPwm} with a slightly different initial guess:
%t, Theorem \ref{IVPwtheorem} will hold true for many other initial approximation in $\mathcal{C}_0^{1,2}([0,T)\times D; \mathbb{R})\cap \mathcal{C}([0,T]\times D;\mathbb{R}) $. As an alternative, the iterative approximation scene can be amended to generate the upper and lower bounds which also converges to the unknown true solution.
%\textcolor{red}{discuss more}
%To this end, we define the sequence $\{u_m\}_{m\in\mathbb{N}_0}$ to start with 
\begin{equation}\label{u0pro}
u_0(t,{\bf x},i):= \mathbb{E}_0^{t,{\bf x},i}\left[\Theta_{t,T}^i\Lambda^i_{t,T}g(X_T,i) -\int_t^T \Theta_{t,s}^i\Lambda^i_{t,s}\phi(s,X_s,i)ds\right],\quad (t,{\bf x},i)\in [0,T]\times D\times \mathcal{M},
\end{equation}
which uniquely solves the initial value problem:
\begin{equation}\label{u0ivp}
\begin{dcases}
\frac{\partial}{\partial t}u_0(t,{\bf x},i)+\mathcal{L}_{i}u_0( t,{\bf x},i)  =  (r(t,{\bf x},i)-q_{ii}({\bf x}))u_0(t,{\bf x},i)+\phi(t,{\bf x},i),
&  (t,{\bf x},i) \in [0,T)\times D  \times \mathcal{M},\\
 u_0(T,{\bf x},i)=g({\bf x},i),&  ({\bf x},i) \in  D  \times \mathcal{M},
\end{dcases}
\end{equation}
under Assumption \ref{Ai}, just as so is \eqref{def of w0} to \eqref{IVPw0pde}, with the only difference being the additional multiple $\Lambda_{t,\cdot}^i$ inside the expectation \eqref{u0pro} and the corresponding potential $-q_{ii}({\bf x})$ in \eqref{u0ivp}.
Then, as mentioned earlier, in exactly the same spirit as the recursion \eqref{IVPwm}, we construct the recursive sequence $\{u_m\}_{m\in\mathbb{N}}$ through either the probabilistic representation:
\begin{equation}\label{umpro0}
u_m(t,{\bf x},i)=\mathbb{E}_0^{t,{\bf x},i}\left[\Theta_{t,T}^i\Lambda^i_{t,T}g(X_T,i) -\int_t^T \Theta_{t,s}^i\Lambda^i_{t,s} \left(\phi(s,X_s,i) -\sum_{j\in\mathcal{M}\setminus\{i\}}q_{ij}(X_s)u_{m-1}(s,X_s,j) \right)ds\right],
\end{equation}
or the initial value problem:
\[
\begin{dcases}
 \frac{\partial}{\partial t} u_m(t,{\bf x},i) +\mathcal{L}_i u_m(t,{\bf x},i) = (r(t,{\bf x},i)-q_{ii}({\bf x}))u_m(t,{\bf x},i)+\phi(t,{\bf x},i)-\sum_{j\in\mathcal{M}\setminus\{i\}}q_{ij}({\bf x})u_{m-1}(t,{\bf x},j),
& (t,{\bf x},i) \in [0,T)\times D \times \mathcal{M}, \\
u_m(T,{\bf x},i) =g({\bf x},i),
& ({\bf x},i) \in D \times \mathcal{M}.
\end{dcases}
\]
As mentioned earlier, from a complexity point of view, this iterative scheme is equivalent to that of Section \ref{section iterative weak approximation}, in the sense that the only difference is the presence of $\Lambda_{t,\cdot}^i$ in \eqref{u0pro}, which costs effectively none.
The motivation behind this alternative scheme is Theorem \ref{IVPpro} {\bf (c)}, that is, it may offer an additional feature of monotonic convergence under suitable conditions.
Notably, the additional condition here is readily verifiable, as it is on the initial condition and the heat source. 

\begin{theorem} \label{IVPpro}
{\bf (a)} The functions $w_m$, $u_m$ and $u_{m-1}$ satisfy the system:
\begin{align}\label{wmpro}
w_m (t,{\bf x},i) &= u_{m-1} (t,{\bf x},i) + \mathbb{E}_Q^{t,{\bf x},i} \left[ \mathbbm{1}(\tau_t^{(m)}\leq T)\Theta_{t,\tau_t^{(m)}}w_0(\tau_t^{(m)},X_{\tau_t^{(m)}},R_{\tau_t^{(m)}}) \right],\\
u_m (t,{\bf x},i) &= u_{m-1} (t,{\bf x},i) + \mathbb{E}_Q^{t,{\bf x},i} \left[ \mathbbm{1}(\tau_t^{(m)}\leq T)\Theta_{t,\tau_t^{(m)}}u_0(\tau_t^{(m)},X_{\tau_t^{(m)}},R_{\tau_t^{(m)}}) \right].\label{umpro}
\end{align}

\noindent {\bf (b)}
It holds that for $(t,{\bf x},i) \in [0,T)\times D \times \mathcal{M}$, $|u_m(t,{\bf x},i)-v(t,{\bf x},i)|=\mathcal{O}((c^{m+1}(T-t)^{m+1}/(m+1)!)^p)$ as $m\to +\infty$ for any $p\in (0,1)$.
If each of the sequences $\{g(\cdot,i)\}_{i\in \mathcal{M}}$ and $\{\phi(\cdot,\cdot,i)\}_{i\in\mathcal{M}}$ is either non-negative or non-positive, then it holds that $u_m(\cdot,\cdot,i)\to  v(\cdot,\cdot,i)$ locally uniformly on $[0,T]\times D$ for $i\in\mathcal{M}$.

\noindent {\bf (c)}
If, moreover, $\{g(\cdot,i)\}_{i\in \mathcal{M}}$ is non-negative (resp. non-positive) and $\{\phi(\cdot,\cdot,i)\}_{i\in\mathcal{M}}$ is non-positive (resp. non-negative), then the locally uniform convergence is monotonic from below $u_m(\cdot,\cdot,i)\uparrow  v(\cdot,\cdot,i)$ (resp. from above $u_m(\cdot,\cdot,i)\downarrow  v(\cdot,\cdot,i)$) as $m\to +\infty$.
\end{theorem}

\begin{proof}[Proof of Theorem \ref{IVPpro}]\label{proofAp}
{\bf (a)} 
Under the probability measure $\mathbb{P}_Q^{t,{\bf x},i}$, the indicator function $\mathbbm{1}(\tau_t^{(1)}>s)$ is infinitely differentiable for all $s\in [t,\tau_t^{(1)}\land T)$, as $\mathbbm{1}(\tau_t^{(1)}>\cdot)\equiv 1$ throughout.
%Moreover, it is clear that $R_s\equiv i$ for all $s\in [t,\tau_t^{(1)}\land T)$.
Therefore, it holds by the Ito formula that for $s\in [t,\tau_t^{(1)}\land T)$,
\begin{align*}
d\left(\mathbbm{1}(\tau_t^{(1)}>s)\Theta_{t,s}u_0(s,X_s,R_s)\right)
&= d\left(\mathbbm{1}(\tau_t^{(1)}>s)\Theta_{t,s}u_0(s,X_s,i)\right)\\
%&=\mathbbm{1}(\tau_t^{(1)}>s)\Theta_{t,s}\left[\phi(s,X_s,R_s)ds+\sum_{j\in\mathcal{M}\setminus \{R_s\}}q_{R_s j}(X_s)u_0(s,X_s,j)+\langle \nabla u_0(s,X_s,R_s),\sigma(s,X_s,R_s)dW_s\rangle\right]\\
&=\mathbbm{1}(\tau_t^{(1)}>s)\Theta_{t,s}\bigg[\left(\frac{\partial}{\partial s}+\mathcal{L}_{i}-r(s,X_s,i)+q_{ii}(X_s)\right)u_0(s,X_s,i)ds\\
&\qquad \qquad \qquad \qquad \qquad +\langle \nabla u_0(s,X_s,i),\sigma(s,X_s,i)dW_s\rangle\bigg]\\
&=\mathbbm{1}(\tau_t^{(1)}>s)\Theta_{t,s}\left[\phi(s,X_s,i)ds+\langle \nabla u_0(s,X_s,i),\sigma(s,X_s,i)dW_s\rangle\right],
\end{align*}
where the third equality holds by \eqref{u0ivp} and $R_s\equiv i$ for all $s\in [t,\tau_t^{(1)}\land T)$ under $\mathbb{P}_Q^{t,{\bf x},i}$.
In the second equality, no inter-dependency terms $\sum_{j\in \mathcal{M}\setminus \{i\}}$ appear anymore as the Ito formula has been applied to $u_0(s,X_s,i)$, that is, with the third argument being degenerate at the regime $i$.
By integrating from $t$ to $\tau_t^{(1)}\land T$ and taking expectation under $\mathbb{P}_Q^{t,{\bf x},i}$, we obtain
\begin{align*}
&\mathbb{E}_Q^{t,{\bf x},i}\left[\mathbbm{1}(\tau_t^{(1)}>{\tau_t^{(1)}\land T})\Theta_{t,{\tau_t^{(1)}\land T}}u_0({\tau_t^{(1)}\land T},X_{{\tau_t^{(1)}\land T}},i)\right]\\
&\qquad =\mathbb{E}_Q^{t,{\bf x},i}\left[\mathbbm{1}(\tau_t^{(1)}>t)\Theta_{t,t}u_0(t,X_t,i)+\int_t^{{\tau_t^{(1)}\land T}} \mathbbm{1}(\tau_t^{(1)}>s)\Theta_{t,s}\phi(s,X_s,i)ds\right]\\
&\qquad =u_0(t,{\bf x},i)+\mathbb{E}_Q^{t,{\bf x},i}\left[\int_t^{\tau_t^{(1)}\land T} \Theta_{t,s}\phi(s,X_s,i)ds\right],
\end{align*}
which yields
\begin{equation}\label{u0 representation}
u_0(t,{\bf x},i)= \mathbb{E}_Q^{t,{\bf x},i} \left[ \mathbbm{1}(\tau_t^{(1)}>T) \Theta_{t,T}g(X_T,i) - \int_t^{\tau_t^{(1)} \land T} \Theta_{t,s}\phi(s,X_s,i)ds\right],
\end{equation}
due to $\mathbbm{1}(s>s\land T)=\mathbbm{1}(s>T)\mathbbm{1}(s>s\land T)+\mathbbm{1}(s\le T)\mathbbm{1}(s>s\land T)=\mathbbm{1}(s>T)\mathbbm{1}(s>T)+\mathbbm{1}(s\le T)\mathbbm{1}(s>s)=\mathbbm{1}(s>T)$ for all $s\in [t,+\infty)$.
Note that we have left $\Theta_{t,T}$ unchanged to $\Theta_{t,T}^i$ in \eqref{u0 representation} on purpose for later convenience.

Next, along the same line as the proof of Theorem \ref{IVPwtheorem} {\bf (a)} with \eqref{wmdef}, \eqref{IVPwm} and \eqref{IVPwmpde}, we obtain the representation:
\begin{equation}\label{um representation}
 u_m(t,{\bf x},i)=\mathbb{E}_Q^{t,{\bf x},i}\left[\mathbbm{1}(\tau_t^{(m)}>T) \Theta_{t,T}g(X_T,R_T)   - \int_t^{\tau_t^{(m)} \land T} \Theta_{t,s}\phi(s,X_s,R_s)ds + \mathbbm{1}(\tau_t^{(m)}\leq T)\Theta_{t,\tau_t^{(m)}}u_0(\tau_t^{(m)},X_{\tau_t^{(m)}},R_{\tau_t^{(m)}})\right],
\end{equation}
for all $m\in\mathbb{N}$.
The third term in \eqref{um representation} can be decomposed as follows:
\begin{align}
&\mathbb{E}_Q^{t,{\bf x},i}\left[\mathbbm{1}(\tau_t^{(m)}\leq T)\Theta_{t,\tau_t^{(m)}}u_0(\tau_t^{(m)},X_{\tau_t^{(m)}},R_{\tau_t^{(m)}})\right]\nonumber\\
&\qquad =\mathbb{E}_Q^{t,{\bf x},i}\left[\mathbbm{1}(\tau_t^{(m)}\leq T)\Theta_{t,\tau_t^{(m)}}\mathbb{E}_Q^{\tau_t^{(m)}, X_{\tau_t^{(m)}},R_{\tau_t^{(m)}}}\left[\mathbbm{1}(\tau_{\tau_t^{(m)}}^{(1)}>T)\Theta_{\tau_t^{(m)},T}g(X_T,R_T)-\int_{\tau_t^{(m)}}^{\tau_{\tau_t^{(m)}}^{(1)} \land T} \Theta_{\tau_t^{(m)},s}\phi(s,X_s,R_s)ds \right]\right]\nonumber\\
&\qquad =\mathbb{E}_Q^{t,{\bf x},i}\left[\mathbbm{1}(\tau_t^{(m)}\leq T)\Theta_{t,\tau_t^{(m)}}\mathbb{E}_Q^{t,{\bf x},i}\left[\mathbbm{1}(\tau_t^{(m+1)}>T)\Theta_{\tau_t^{(m)},T}g(X_T,R_T)-\int_{\tau_t^{(m)}}^{\tau_t^{(m+1)} \land T} \Theta_{\tau_t^{(m)},s}\phi(s,X_s,R_s)ds \Bigg|\,\mathcal{F}_{\tau_t^{(m)}}\right]\right]\nonumber\\
&\qquad =\mathbb{E}_Q^{t,{\bf x},i}\left[\mathbb{E}_Q^{t,{\bf x},i}\left[\mathbbm{1}(\tau_t^{(m)}\leq T,\,\tau_t^{(m+1)}>T)\Theta_{t,T}g(X_T,R_T)-\mathbbm{1}(\tau_t^{(m)}\leq T)\int_{\tau_t^{(m)}}^{\tau_t^{(m+1)} \land T} \Theta_{t,s}\phi(s,X_s,R_s)ds \Bigg|\,\mathcal{F}_{\tau_t^{(m)}}\right]\right]\nonumber\\
&\qquad =\mathbb{E}_Q^{t,{\bf x},i}\left[\mathbbm{1}(\tau_t^{(m)}\leq T,\,\tau_t^{(m+1)}>T)\Theta_{t,T}g(X_T,R_T)-\mathbbm{1}(\tau_t^{(m)}\leq T)\int_{\tau_t^{(m)}}^{\tau_t^{(m+1)} \land T} \Theta_{t,s}\phi(s,X_s,R_s)ds \right],\label{second term in detail}
\end{align}
where the first equality holds by \eqref{u0 representation}, the second by the strong Markov property, the third by $\mathcal{F}_{\tau_t^{(m)}}$-measurability of the random variable $\mathbbm{1}(\tau_t^{(m)}\leq T)\Theta_{t,\tau_t^{(m)}}$, and the last by the tower property.
Combining \eqref{um representation} and \eqref{second term in detail} yields
\begin{equation}\label{um m+1}
 u_m(t,{\bf x},i)=\mathbb{E}_Q^{t,{\bf x},i}\left[\mathbbm{1}(\tau_t^{(m+1)}>T) \Theta_{t,T}g(X_T,R_T)-\int_t^{\tau_t^{(m+1)} \land T} \Theta_{t,s}\phi(s,X_s,R_s)ds\right], \quad (t,{\bf x},i)\in [0,T]\times D\times \mathcal{M},\quad m\in\mathbb{N},
\end{equation}
due to $\mathbbm{1}(\tau_t^{(m)}>T)+\mathbbm{1}(\tau_t^{(m)}\leq T,\,\tau_t^{(m+1)}>T)=\mathbbm{1}(\tau_t^{(m+1)}>T)$ and $\int_t^{\tau_t^{(m)} \land T}+\mathbbm{1}(\tau_t^{(m)}\leq T)\int_{\tau_t^{(m)}}^{\tau_t^{(m+1)} \land T}=\int_t^{\tau_t^{(m+1)}\land T}$.
Hence, the recursions \eqref{wmpro} and \eqref{umpro} can be derived by rewriting \eqref{wmdef} and \eqref{um representation} with $u_{m-1}$ on the basis of the representation \eqref{um m+1}.

\noindent
{\bf (b)} In a similar manner to the proof of Theorem \ref{IVPwtheorem}, the point-wise convergence and the decay rate holds by the dominated convergence theorem with the aid of \cite[Proposition 2.3]{YZ1}, due to $\tau_t^{(m)}$ diverges in $m$ by Assumption \ref{qproperty} {\bf (c)}.

\noindent 
{\bf (c)}
To ease the notation, we write $v(t,{\bf x},i)-u_m(t,{\bf x},i)=v_{m,1}(t,{\bf x},i)-v_{m,2}(t,{\bf x},i)$, where
\[
 v_{m,1}(t,{\bf x},i):=\mathbb{E}_Q^{t,{\bf x},i}\left[\mathbbm{1}(\tau_t^{(m+1)}\le T) \Theta_{t,T}g(X_T,R_T)\right],\quad
 v_{m,2}(t,{\bf x},i):=\mathbb{E}_Q^{t,{\bf x},i}\left[\int_{\tau_t^{(m+1)} \land T}^T \Theta_{t,s}\phi(s,X_s,R_s)ds\right],
\]
for $(t,{\bf x},i)\in [0,T]\times D \times \mathcal{M}$ and $m\in\mathbb{N}$.
If each of the sequences $\{g(\cdot,i)\}_{i\in \mathcal{M}}$ and $\{\phi(\cdot,\cdot,i)\}_{i\in\mathcal{M}}$ is either non-negative or non-positive, then $\{v_{m,1}(\cdot,\cdot,i)\}_{m\in\mathbb{N}}$ and $\{v_{m,2}(\cdot,\cdot,i)\}_{m\in\mathbb{N}}$ are monotonic sequences of real-valued continuous functions for each $i\in\mathcal{M}$.
Hence, $\{u_m\}_{m\in\mathbb{N}}$ converges to $v$ locally uniformly by the Dini theorem.
If, moreover, the two functions $g$ and $\phi$ have opposite signs (where zero is allowed), then the sequence $\{v_{m,1}(\cdot,\cdot,i)-v_{m,2}(\cdot,\cdot,i)\}_{m\in\mathbb{N}}$ on the whole is monotonic.
Finally, in comparison of \eqref{defv} and \eqref{um m+1}, it is straightforward to see whether from above or below.
\end{proof}

\subsection{Hard bounding functions}\label{section hard bounding functions IVP}

In the context of iterative methods, it is almost as imperative as the method per se to equip it with a valid stopping criterion.
In the literature, there exist a variety of such criteria, whereas the majority is largely based on some error from one iteration to the other, that is, the iteration terminates when the error is found below a predetermined tolerance.
The final product is expected to be very close to the unknown target, provided that a convergence towards the target is guaranteed, whereas it is still an approximation, for instance, with no ability of identifying whether it sits above or below the unknown target in any way.  
In our case, we have derived a convergence rate in Theorem \ref{IVPwtheorem} {\bf (b)}, whereas it does not seem very useful for setting a practical stopping criterion either, as it is not an exact rate as well as due to a lack of computable constant multiples. 

Here, we instead employ the approach taken in \cite{rayK} to construct hard upper and lower bounding functions, not approximations, for the unknown target function \eqref{defv} on the basis of the obtained sequences of approximations $\{w_m\}_{m\in\mathbb{N}_0}$ and $\{u_m\}_{m\in\mathbb{N}_0}$.
Remarkably, a paired sequence of the computable hard upper and lower bounding functions converge towards each other, that is, the target function \eqref{defv} will be fully specified eventually without an exception.
Given that those upper and lower bounding functions are easily computable, one may set up very reliable stopping criteria on the basis of the distance between upper and lower bounds, rather than the distance between two consecutive approximations.

We prepare some notation.
For each $m\in\mathbb{N}_0$, we define the following two functions on $[0,T]\times D \times \mathcal{M}$:
\begin{equation}\label{def of Mr}
M_r(t,{\bf x},i):=\sum_{j\in\mathcal{M}}q_{ij}({\bf x})\mathbb{E}_0^{t,{\bf x},j}\left[\int^T_t \Theta_{t,s}^j ds\right],
%=\sum_{j\in\mathcal{M}}q_{ij}({\bf x})\mathbb{E}_0^{t,{\bf x},j}\left[\int^T_t \exp\left[-\int_t^s r(u,X_u,j)du \right] ds\right],
\end{equation}
and, for an array $\{f_m(\cdot,\cdot,j)\}_{j\in\mathcal{M},\,m\in\mathbb{N}_0}$ of functions on $[0,T]\times D$,
\begin{equation}\label{bound functions}
N_m(t,{\bf x},i;\,f):=
\begin{dcases}
\sum_{j\in\mathcal{M}}q_{ij}({\bf x})f_0(t,{\bf x},j), &\text{if } m=0,\\
\sum_{j\in\mathcal{M}\setminus\{i\}}q_{ij}({\bf x})\left(f_m(t,{\bf x},j)-f_{m-1}(t,{\bf x},j)\right), &\text{if } m\in \mathbb{N}.
\end{dcases}
\end{equation}
%where the last equality holds because no regime switching can occur under $\mathbb{P}_0^{t,{\bf x},j}$.
Moreover, we define $M_r^U$, $M_r^L$, $\{N_m^U(f)\}_{m\in\mathbb{N}_0}$ and $\{N_m^L(f)\}_{m\in\mathbb{N}_0}$ as the essential supremum and infimum of $M_r$ and $\{N_m(\cdot,\cdot,\cdot;\,f)\}_{m\in\mathbb{N}_0}$, as follows:
\begin{gather}
M_r^U:=\esssup_{(t,{\bf x},i)\in[0,T]\times D \times\mathcal{M}}\left(M_r(t,{\bf x},i)\right)_+, \quad 
M_r^L:=\essinf_{(t,{\bf x},i)\in[0,T]\times D \times\mathcal{M}}\left(M_r(t,{\bf x},i)\right)_-,\label{def of MrU MrL}\\
N_m^U(f):=\esssup_{(t,{\bf x},i)\in[0,T]\times D \times\mathcal{M}}\left(N_m(t,{\bf x},i;\,f)\right)_+, \quad 
N_m^L(f):=\essinf_{(t,{\bf x},i)\in[0,T]\times D \times\mathcal{M}}\left(N_m(t,{\bf x},i;\,f)\right)_-,\label{def of NmU NmL}
\end{gather}
and define $\{U_m(\cdot,\cdot,\cdot;\,f)\}_{m\in\mathbb{N}}$ and $\{L_m(\cdot,\cdot,\cdot;\,f)\}_{m\in\mathbb{N}}$ by
\begin{equation}\label{LUm}
U_m(t,{\bf x},i;\,f):=f_m(t,{\bf x},i)+\frac{N_m^U(f)}{1-M_r^U}\mathbb{E}_{0}^{t,{\bf x},i} \left[\int_{t}^{T}\Theta_{t,s}^i ds\right],\quad 
L_m(t,{\bf x},i;\,f):=f_m(t,{\bf x},i)+\frac{N_m^L(f)}{1-M_r^L}\mathbb{E}_{0}^{t,{\bf x},i} \left[\int_{t}^{T}\Theta_{t,s}^i ds\right].
\end{equation}
We are now in a position to give the main result of this section.

\begin{theorem}\label{IVPB}
Let $f$ be either $\{w_m\}_{m\in\mathbb{N}_0}$ or $\{u_m\}_{m\in\mathbb{N}_0}$.

\noindent {\bf (a)}
It holds that $|N_m(t,{\bf x},i;\,f)| \to  0$ as $m\to +\infty$ for all $(t,{\bf x},i)\in[0,T]\times D \times\mathcal{M}$.

\noindent {\bf (b)} 
If $M_r^U<1$ and $|N_m^U(f)|+|N_m^L(f)|<+\infty$ for all $m\in\mathbb{N}_0$, then it holds that
\begin{equation}\label{boundswm0}
L_m(t,{\bf x},i;\,f)\le v(t,{\bf x},i) \le U_m(t,{\bf x},i;\,f),\quad (t,{\bf x},i)\in [0,T]\times D\times \mathcal{M}.
\end{equation}

\noindent {\bf (c)}
If, moreover, the conditions of Theorem \ref{IVPpro} {\bf (c)} holds true, then we have $u_m(t,{\bf x},i)\le v(t,{\bf x},i)\le U_m(t,{\bf x},i;\,u)$ (resp. $L_m(t,{\bf x},i;u)\le v(t,{\bf x},i)\le u_m(t,{\bf x},i)$), where both upper and lower bounds tend to $v$ locally uniformly on $[0,T]\times D$ for $i\in \mathcal{M}$ as $m\to +\infty$.
\end{theorem}

Some trivial yet interesting properties follow immediately from \eqref{LUm}.
For instance, if $w_m$ converges to $v$ uniformly (Theorem \ref{IVPwtheorem} {\bf (b)}), then both $N^U_m(w)$ and $N^L_m(w)$ tend to vanish due to the structure \eqref{bound functions}.
Hence, with the aid of the boundedness $\Theta_{t,s}\in [0,1]$, the bounding functions $L_m$ and $U_m$ both converge uniformly to the target solution $v$ as $m\to +\infty$.

\begin{proof}[Proof of Theorem \ref{IVPB}]
To ease the notation, we employ a set of differential operators $\{\mathcal{A}_{i}: i\in\mathcal{M}\}$ on a collection of smooth functions $\{f(\cdot,\cdot,k):\,k\in\mathcal{M}\}$, defined by
\begin{equation}\label{operator Ai}
\mathcal{A}_{i} f(t,{\bf x},k) := \mathcal{L}_i  f(t,{\bf x},k)-r(t,{\bf x},k)f(t,{\bf x},k) + \sum_{j\in\mathcal{M}}q_{ij}({\bf x}) f(t,{\bf x},j). %\quad (t,{\bf x})\in[0,T)\times\mathbb{R}^n.
\end{equation}
Although the operator $\mathcal{A}_{i}$ can be applied to any member of the collection, we focus on the case $i=k$, because the other cases are of no value in our context.

We prove the results only for $\{w_m\}_{m\in\mathbb{N}_0}$, as that for $\{u_m\}_{m\in\mathbb{N}_0}$ is rather repetitious.
We get the lower bound for $v(t,{\bf x},i)$ as follows:
\begin{align*}
v(t,{\bf x},i)&=\mathbb{E}_Q^{t,{\bf x},i}\left[\Theta_{t,T}g(X_T,R_T) - \int_t^T \Theta_{t,s} \phi(s,X_s,R_s) ds \right] \\
& \ge \mathbb{E}_Q^{t,{\bf x},i}\left[\Theta_{t,T}w_0(T,X_T,R_T)-\int_t^T \Theta_{t,s}\left(\frac{\partial}{\partial s}w_0(s,X_s,R_s)+ \mathcal{A}_{R_s}w_0(s,X_s,R_s) - N_0^L(w)\right)ds\right]\\
&=w_0(t,{\bf x},i) + N_0^L(w) \mathbb{E}_Q^{t,{\bf x},i}\left[\int_t^T \Theta_{t,s}ds\right],
\end{align*}
where we have applied $g(X_T,R_T)=w_0(T,X_T,R_T)$ and for $s\in[t,T)$,
\begin{align*}
\phi(s,X_s,R_s)
&=\frac{\partial}{\partial s}w_0(s,X_s,R_s)+ \mathcal{A}_{R_s}w_0(s,X_s,R_s)- \sum_{j\in\mathcal{M}}q_{ R_s j}(X_s)w_0(s,X_s,j)\\
&\leq \frac{\partial}{\partial s}w_0(s,X_s,R_s)+ \mathcal{A}_{R_s}w_0(s,X_s,R_s) - N_0^L(w).
\end{align*}
due to \eqref{IVPw0pde}.
Similarly, we obtain the upper bound for $v(t,{\bf x},i)$, which yields the inequalities:
\[
N_0^L(w) \mathbb{E}_Q^{t,{\bf x},i}\left[\int_t^T \Theta_{t,s}ds\right]
\leq v(t,{\bf x},i)-w_0(t,{\bf x},i)
\leq N_0^U(w) \mathbb{E}_Q^{t,{\bf x},i}\left[\int_t^T \Theta_{t,s}ds\right].
\]
By setting $g\equiv 0$ and $\phi\equiv-1$ here, we obtain the following inequalities:
\[
M_r^L\mathbb{E}_Q^{t,{\bf x},i}\left[\int_t^T \Theta_{t,s}ds\right]
\leq \mathbb{E}_Q^{t,{\bf x},i}\left[\int_t^T \Theta_{t,s}ds\right] - \mathbb{E}_0^{t,{\bf x},i}\left[\int_t^T \Theta_{t,s}^ids\right]
\leq M_r^U\mathbb{E}_Q^{t,{\bf x},i}\left[\int_t^T \Theta_{t,s}ds\right].
\]
By rearranging the inequalities above under the condition $M_r^U<1$, we obtain the inequalities:
\begin{equation}\label{boundsr}
\frac{1}{1-M_r^L}\mathbb{E}_0^{t,{\bf x},i}\left[\int_t^T \Theta_{t,s}^i ds\right]
\leq \mathbb{E}_Q^{t,{\bf x},i}\left[\int_t^T \Theta_{t,s}ds\right]
\leq \frac{1}{1-M_r^U}\mathbb{E}_0^{t,{\bf x},i}\left[\int_t^T \Theta_{t,s}^ids\right].
\end{equation}
Combining the two results yields the desired inequalities \eqref{boundswm0} for $m=0$. 
The case for $m\in\mathbb{N}$ proceeds in a similar manner, where the lower bound can be obtain for $v(t,{\bf x},i)$ as follows:
\begin{align*}
v(t,{\bf x},i)&=\mathbb{E}_Q^{t,{\bf x},i}\left[\Theta_{t,T}g(X_T,R_T) - \int_t^T \Theta_{t,s} \phi(s,X_s,R_s) ds \right] \\
&\ge \mathbb{E}_Q^{t,{\bf x},i}\left[\Theta_{t,T}w_m(T,X_T,R_T)-\int_t^T \Theta_{t,s}\left(\frac{\partial}{\partial s}w_m(s,X_s,R_s)+ \mathcal{A}_{R_s}w_m(s,X_s,R_s) - N_m^L(w)\right)ds\right]\\
&=w_m(t,{\bf x},i) + N_m^L(w) \mathbb{E}_Q^{t,{\bf x},i}\left[\int_t^T \Theta_{t,s}ds\right],
\end{align*}
where we have applied $g(X_T,R_T)=w_m(T,X_T,R_T)$ and for $s\in[t,T)$,
\begin{align*}
\phi(s,X_s,R_s) &= \frac{\partial}{\partial s}w_m(s,X_s,R_s)+ \mathcal{L}_{R_s}w_m(s,X_s,R_s) -(r(s,X_s,R_s)-q_{R_s R_s}(X_s))w_m(s,X_s,R_s) + \sum_{j\in\mathcal{M}\setminus\{R_s\}}q_{R_s j}(X_s)w_{m-1}(s,X_s,R_s) \\
%&= \frac{\partial}{\partial s}w_m(s,X_s,R_s)+ \mathcal{A}_{R_s,1}w_m(s,X_s,R_s) -r(s,X_s,R_s)w_m(s,X_s,R_s) -\sum_{j\in\mathcal{M}\setminus\{R_s\}}q_{R_s j}(X_s)\left(w_m(s,X_s,R_s)-w_{m-1}(s,X_s,R_s)\right) \\
&\leq  \frac{\partial}{\partial s}w_m(s,X_s,R_s)+ \mathcal{A}_{R_s}w_m(s,X_s,R_s)- N_m^L(w),
\end{align*}
due to \eqref{IVPwmpde}.
This yields the inequalities
\begin{equation}\label{boundswmQ}
N_m^L(w) \mathbb{E}_Q^{t,{\bf x},i}\left[\int_t^T \Theta_{t,s}ds\right]
\leq v(t,{\bf x},i)-w_m(t,{\bf x},i)
\leq N_m^U(w) \mathbb{E}_Q^{t,{\bf x},i}\left[\int_t^T \Theta_{t,s}ds\right].
\end{equation}
Combining \eqref{boundsr} and \eqref{boundswmQ} yields the desired inequalities \eqref{boundswm0} for all $m\in\mathbb{N}_0$.
Since $\mathcal{M}$ is a finite set and $w_m$ converges to $v$ for all $(t,{\bf x},i)\in[0,T]\times D\times\mathcal{M}$, we have
\begin{align*}
\left|N_m(t,{\bf x},i;\,w)\right|&\leq\sum_{j\in\mathcal{M}\setminus\{i\}}q_{ij}({\bf x})\left|w_m(t,{\bf x},j)-w_{m-1}(t,{\bf x},j)\right|\\
&\leq \sum_{j\in\mathcal{M}\setminus\{i\}}q_{ij}({\bf x})\left|w_m(t,{\bf x},j)-v(t,{\bf x},j)\right|
+\sum_{j\in\mathcal{M}\setminus\{i\}}q_{ij}({\bf x})\left|w_{m-1}(t,{\bf x},j)-v(t,{\bf x},j)\right|\to  0,
\end{align*}
as $m\to +\infty$, which concludes the proof.
\end{proof}

It is worth noting that the computation of the bounding functions \eqref{LUm} simplifies under some practical conditions. 
For instance, if the function $r(t,\cdot,i)$ is independent of the state variable ${\bf x}$, then the bounding functions \eqref{LUm} reduce to
\begin{equation}\label{LUmsimple0}
U_m(t,{\bf x},i;\,f)=f_m(t,{\bf x},i)+\frac{N_m^U(f)}{1-M_r^U} \int_t^T\Theta_{t,s}^i ds,\quad
L_m(t,{\bf x},i;\,f)=f_m(t,{\bf x},i)+\frac{N_m^L(f)}{1-M_r^L} \int_t^T\Theta_{t,s}^i ds,
\end{equation} 
which skips one expectation, as $\Theta_{t,s}$ has no randomness anymore.
If moreover the function $r(t,\cdot,\cdot)$ is independent of the regime variable, then we have $M_r\equiv 0$ due to Assumption \ref{qproperty} {\bf (b)}, and thus the bounding functions \eqref{LUm} reduce to
\begin{equation}\label{LUmsimple}
U_m(t,{\bf x},i;\,f)=f_m(t,{\bf x},i)+N_m^U(f) \int_t^T\Theta_{t,s}^i ds,\quad 
L_m(t,{\bf x},i;\,f)=f_m(t,{\bf x},i)+N_m^L(f) \int_t^T\Theta_{t,s}^i ds,
\end{equation} 
which allow one to skip the computation for \eqref{def of Mr} and \eqref{def of MrU MrL} at all.

\section{Initial boundary value problems}\label{section initial boundary value problems}

In this section, we briefly discuss the proposed iterative weak approximation and hard bounding functions in the case where there exist attainable boundaries at which the regime-switching diffusion \eqref{XSDE} ends its life. 
Recall that the first hitting time to such boundaries is defined by $\eta_t:=\inf\{s\geq t:\,X_s \not\in D\}$ in \eqref{etadef}, with a bounded open domain $D(\subset \mathbb{R}^n)$.
In order to avoid overloading the paper, we keep this section as concise as possible, for instance, by omitting proofs, because the problem setting and the results are considerably similar to those in Section \ref{section initial value problems}, or even less sensitive to various factors in proving the results, thanks to the compactness of the domain.

The target solution that we consider in this section is written in the terms of mathematical expectation:
\begin{equation}\label{IBdefv}
v(t,{\bf x},i):= \mathbb{E}_Q^{t,{\bf x},i} \left[ \mathbbm{1}(\eta_t\geq T)\Theta_{t,T}  g(X_T, R_T)
+\mathbbm{1}(\eta_t <T) \Theta_{t,\eta_t}\Psi(\eta_t, X_{\eta_t},R_{\eta_t})
 - \int_t^{\eta_t \land T} \Theta_{t,s}\phi(s,X_s,R_s)ds \right],
\end{equation}
for all $(t,{\bf x},i) \in  [0,T]\times\overline{D} \times \mathcal{M}$, which is a unique solution to the initial boundary value problem:
\begin{equation}\label{IBpdev}
\begin{dcases}
\frac{\partial}{\partial t}v(t,{\bf x},i)+\mathcal{L}_{i}v( t,{\bf x},i) =  (r(t,{\bf x},i)-q_{ii}({\bf x}))v(t,{\bf x},i)+\phi(t,{\bf x},i)-\sum_{j\in\mathcal{M}\setminus\{i\}}q_{ij}({\bf x})v(t,{\bf x},j),
&  (t,{\bf x},i) \in [0,T)\times D  \times \mathcal{M};\\
 v(T,{\bf x},i)=g({\bf x},i),
&  ({\bf x},i) \in  D  \times \mathcal{M};\\
 v(t,{\bf x},i) = \Psi(t,{\bf x},i),
& (t,{\bf x},i) \in [0,T] \times \partial D \times \mathcal{M},
\end{dcases}
\end{equation}
under the following conditions on the data \cite{BYZ, YZ1}.

\begin{assumption}\label{Aib}
There exists an $\alpha\in(0,1)$, such that for each $i\in\mathcal{M}$,

\noindent 
{\bf (a)}  $r(\cdot,\cdot,i)$ is $\alpha$-H\"older continuous on $[0,T]\times\overline{D}$,

\noindent
{\bf (b)}   $\phi(\cdot,{\bf x},i)$ is $\alpha$-H\"older continuous on $[0,T]\times\overline{D}$,

\noindent 
{\bf (c)}  $g(\cdot,i)$ is continuous on $\overline{D}$, and

\noindent
{\bf (d)} $\Psi(\cdot,\cdot,i)$ is continuous on $\partial D \times [0,T]$, with $g({\bf x},i)=\Psi(T,{\bf x},i)$ for all ${\bf x}\in \partial D$.
\end{assumption}

%\subsection{Iterative weak approximation}\label{section IB iterative}

In a similar spirit to Section \ref{section initial value problems}, we construct a convergent iteration, starting with the following function on $[0,T] \times \overline{D}\times \mathcal{M}$:
%iteration method to approximate the expectation $v(t,{\bf x},i)$ by restricting the maximum number of switchings occur in the time interval $(t,\eta_t\land T]$. First of all, we define the initial approximation  $w_0:[0,T] \times D\times \mathcal{M}\to  \mathbb{R}$ as an expectation under $\mathbb{P}_0^{t,{\bf x},i}$,
\begin{equation}\label{def of w0 IVBP}
w_0(t,{\bf x},i):= \mathbb{E}_0^{t,{\bf x},i} \left[ \mathbbm{1}(\eta_t\geq T)\Theta_{t,T}^i  g(X_T, R_T)
+\mathbbm{1}(\eta_t <T) \Theta_{t,\eta_t}^i\Psi(\eta_t, X_{\eta_t},R_{\eta_t})
 - \int_t^{\eta_t \land T} \Theta^i_{t,s}\phi(s,X_s,R_s)ds \right],
\end{equation}
under the probability measure $\mathbb{P}_0^{t,{\bf x},i}$, that is, no switching is possible with $R_s\equiv i$ throughout its life.
This function uniquely solves the initial boundary value problem:
\begin{equation}\label{IBw0pde}
\begin{dcases}
\frac{\partial}{\partial t}w_0(t,{\bf x},i)+\mathcal{L}_{i}w_0( t,{\bf x},i) =  r(t,{\bf x},i)w_0(t,{\bf x},i)+\phi(t,{\bf x},i),
&  (t,{\bf x},i) \in [0,T)\times D  \times \mathcal{M},\\
 w_0(T,{\bf x},i)=g({\bf x},i),
&  ({\bf x},i) \in  D  \times \mathcal{M},\\
w_0(t,{\bf x},i) = \Psi(t,{\bf x},i),
& (t,{\bf x},i) \in [0,T] \times \partial D \times \mathcal{M}.
\end{dcases}
\end{equation}
Then, with $w_0$ as an initial guess, we construct the sequence $\{w_m\}_{m\in\mathbb{N}}$ of continuous functions by recursion:
\begin{multline}
w_m(t,{\bf x},i)=\mathbb{E}_0^{t,{\bf x},i}\Bigg[ 
\mathbbm{1}(T\leq \eta_t) \Theta_{t,T}^i\Lambda^i_{t,T}g(X_T,i) +
\mathbbm{1}(T>\eta_t) \Theta_{t,\eta_t}^i \Lambda^i_{t,\eta_t} \Psi(\eta_t,X_{\eta_t},R_{\eta_t})\\
-\int_t^{\eta_t\land T} \Theta_{t,s}^i\Lambda^i_{t,s}\left(\phi(s,X_s,i)-\sum_{j\in\mathcal{M}\setminus\{i\}}q_{ij}(X_s)w_{m-1}(s,X_s,j)\right)ds\Bigg],\label{IBPwm without switching}
\end{multline}
which uniquely solves the initial boundary value problem:
\[%begin{equation}\label{IBwmpde}
\begin{dcases}
 \frac{\partial}{\partial t} w_m(t,{\bf x},i) +\mathcal{L}_i w_m(t,{\bf x},i) = (r(t,{\bf x},i)-q_{ii}({\bf x}))w_m(t,{\bf x},i)+\phi(t,{\bf x},i)-\sum_{j\in\mathcal{M}\setminus\{i\}}q_{ij}({\bf x})w_{m-1}(t,{\bf x},j),
& (t,{\bf x},i) \in [0,T)\times D \times \mathcal{M},\\
w_m(T,{\bf x},i) =g({\bf x},i),
& ({\bf x},i) \in D \times \mathcal{M},\\
w_m(t,{\bf x},i) = \Psi(t,{\bf x},i),
&(t,{\bf x},i) \in [0,T] \times \partial D \times \mathcal{M}.
\end{dcases}
\]%end{equation}
provided that the term $\phi(t,{\bf x},i)-\sum_{j\in\mathcal{M}\setminus\{i\}}q_{ij}({\bf x})w_{m-1}(t,{\bf x},j)$ is smooth enough (Assumptions \ref{qproperty} and \ref{Aib}) to be treated as a heat source on the whole.

\begin{theorem}\label{IBwtheorem}
Suppose that the functions $w_m$ defined in \eqref{def of w0 IVBP} and \eqref{IBwmdef} are in $\mathcal{C}_0^{1,2}([0,T)\times D;\mathbb{R})\cap \mathcal{C}([0,T]\times \overline{D};\mathbb{R})$  for all $m\in\mathbb{N}_0$ and $i\in\mathcal{M}$.

\noindent {\bf (a)} It holds that for $m\in\mathbb{N}$ and $(t,{\bf x},i) \in [0,T]\times \overline{D} \times \mathcal{M}$,
\begin{equation}\label{IBwmdef}
w_m(t,{\bf x},i):=\mathbb{E}_Q^{t,{\bf x},i}\left[
\Theta_{t,\eta_t\land T\land \tau_t^{(m)}}w_0(\eta_t\land T\land \tau_t^{(m)},X_{\eta_t\land T\land \tau_t^{(m)}},R_{\eta_t\land T\land \tau_t^{(m)}})
-\int_t^{\eta_t \land T\land \tau_t^{(m)}} \Theta_{t,s} \phi(s,X_s,R_s) ds\right].
\end{equation}

\noindent {\bf (b)} It holds that for $(t,{\bf x},i) \in [0,T]\times D \times\mathcal{M}$, $|w_m(t,{\bf x}, i)-v(t,{\bf x},i)|=\mathcal{O}((c^m(T-t)^m/m!)^p)$ as $m\to +\infty$ for any $p\in (0,1)$, with $c:=\sup_{{\bf x}\in D,j\ne k}q_{jk}({\bf x})$.
Moreover, if the functions $\{w_0(\cdot,\cdot,i)\}_{i\in\mathcal{M}}$ are all either non-negative or non-positive on $[0,T]\times D\times\mathcal{M}$, then we have $w_m(\cdot, \cdot,i)\to v(\cdot,\cdot,i)$ locally uniformly on $[0,T]\times D$.
\end{theorem}

%We omit the proof to avoid overloading the paper, as the results can be proved in a similar manner to Theorem \ref{IVPwtheorem}. 
As before, each switching-free iterate $w_m$ in the form \eqref{IBPwm without switching} represents the expression \eqref{IBwmdef} with switching, where the underlying dynamics terminates if the $m$-th switching occurs before the terminal time or the first exit time.
Also, the representation \eqref{IBwmdef} implies that the sequence $\{w_m\}_{m\in\mathbb{N}}$ tends to the target function \eqref{IBdefv} as desired.

%We obtain the recursive relation of the sequence $\{w_m\}$ in the following lemma.
%\begin{lemma}\label{IBwmlemma1}
%For each $m\in\mathbb{N}$, and for  each $(t,{\bf x},i)\in [0,T]\times\overline{D}\times\mathcal{M}$, the function $w_m$ can be rewritten as:
%\begin{equation}\label{IBwmlemma}
%w_m(t,{\bf x},i)=\mathbb{E}_Q^{t,{\bf x},i}\left[  \Theta_{t,T \land \eta_t\land \tau_t^{(1)}}w_{m-1}(T \land \eta_t\land \tau_t^{(1)},X_{T \land \eta_t\land \tau_t^{(1)}},R_{T \land \eta_t\land \tau_t^{(1)}})
%-\int_t^{\eta_t\land T \land \tau_t^{(1)}} \Theta_{t,s} \phi(s,X_s,R_s) ds \right].
%\end{equation}
%\end{lemma}
%We omit the proof of this lemma, as it is very similar to the proof of Lemma \ref{wmlemma1} by replacing $T$ with $\eta_t\land T$.
%In the following main theorem, we show that we can use the function $w_{m-1}$ to compute $w_m$, and the sequence $\{w_m\}$ converges to $v$.
%Recall the notation $\Theta_{t_1,t_2}^i$ and $\Lambda_{t_1,t_2}^i$ defined in \eqref{def of thetai}.
%For $0\leq t_1\leq t_2\leq T$ and $i\in\mathcal{M}$, we denote
%\[
%\Theta_{t_1,t_2}^i:=\exp\left[-\int_{t_1}^{t_2} r(s,X_s,i) ds\right], \quad
%\Lambda_{t_1,t_2}^i:=\exp\left[\int_{t_1}^{t_2} q_{ii}(X_s)ds\right].
%\]

%\subsection{Monotonically convergent iterative weak approximation}

Here, we propose an alternative iterative scheme on the basis of exactly the same recursion (\ref{IBPwm without switching}) with a slightly different initial guess:
\begin{equation}\label{def of u0 IVBP}
u_0(t,{\bf x},i):= \mathbb{E}_0^{t,{\bf x},i} \left[ \mathbbm{1}(\eta_t\geq T)\Theta_{t,T}^i\Lambda_{t,T}^i  g(X_T, R_T)
+\mathbbm{1}(\eta_t <T) \Theta_{t,\eta_t}^i\Lambda_{t,\eta_t}^i\Psi(\eta_t, X_{\eta_t},R_{\eta_t})
 - \int_t^{\eta_t \land T} \Theta^i_{t,s}\Lambda_{t,T}^i\phi(s,X_s,R_s)ds \right],
\end{equation}
which uniquely solves the initial value problem:
\begin{equation}\label{IBu0pde}
\begin{dcases}
\frac{\partial}{\partial t}u_0(t,{\bf x},i)+\mathcal{L}_{i}u_0( t,{\bf x},i) = ( r(t,{\bf x},i)-q_{ii}({\bf x}))u_0(t,{\bf x},i)+\phi(t,{\bf x},i),
&  (t,{\bf x},i) \in [0,T)\times D  \times \mathcal{M},\\
 u_0(T,{\bf x},i)=g({\bf x},i),
&  ({\bf x},i) \in  D  \times \mathcal{M},\\
u_0(t,{\bf x},i) = \Psi(t,{\bf x},i),
& (t,{\bf x},i) \in [0,T] \times \partial D \times \mathcal{M}.
\end{dcases}
\end{equation}
under Assumption \ref{Aib}, just as so is (\ref{def of w0 IVBP}) and (\ref{IBw0pde}), with the only difference being the additional multiple $\Lambda^i_{t,\cdot}$ inside the expectation (\ref{def of u0 IVBP}) and the corresponding potential $-q_{ii}({\bf x})$ in (\ref{IBu0pde}). Then, as mentioned earlier, in exactly the same spirit as the recursion (\ref{IBPwm without switching}), we construct the recursive sequence $\{u_m\}_{m\in\mathbb{N}}$ through either the probabilistic representation:
\begin{multline}
u_m(t,{\bf x},i)=\mathbb{E}_0^{t,{\bf x},i}\Bigg[ 
\mathbbm{1}(\eta_t \geq T) \Theta_{t,T}^i\Lambda^i_{t,T}g(X_T,i) +
\mathbbm{1}(\eta_t < T) \Theta_{t,\eta_t}^i \Lambda^i_{t,\eta_t} \Psi(\eta_t,X_{\eta_t},R_{\eta_t})\\
-\int_t^{\eta_t\land T} \Theta_{t,s}^i\Lambda^i_{t,s}\left(\phi(s,X_s,i)-\sum_{j\in\mathcal{M}\setminus\{i\}}q_{ij}(X_s)u_{m-1}(s,X_s,j)\right)ds\Bigg],\label{IBPum without switching}
\end{multline}
or the initial boundary value problem:
\begin{equation}
\begin{dcases}
 \frac{\partial}{\partial t} u_m(t,{\bf x},i) +\mathcal{L}_i u_m(t,{\bf x},i) = (r(t,{\bf x},i)-q_{ii}({\bf x}))u_m(t,{\bf x},i)+\phi(t,{\bf x},i)-\sum_{j\in\mathcal{M}\setminus\{i\}}q_{ij}({\bf x})u_{m-1}(t,{\bf x},j),
& (t,{\bf x},i) \in [0,T)\times D \times \mathcal{M},\\
u_m(T,{\bf x},i) =g({\bf x},i),
& ({\bf x},i) \in D \times \mathcal{M},\\
u_m(t,{\bf x},i) = \Psi(t,{\bf x},i),
&(t,{\bf x},i) \in [0,T] \times \partial D \times \mathcal{M}.
\end{dcases}
\end{equation}

As mentioned earlier, from a complexity point of view, this iterative scheme is equivalent to that for the sequence $\{w_m\}_{m\in\mathbb{N}}$, in the sense that the only difference is the presence of $\Lambda_{t,\cdot}^i$ in \eqref{def of u0 IVBP}, which costs effectively none. The motivation behind this alternative scheme is Theorem \ref{IBPpro} {\bf (c)}, that is, it may offer an additional feature of monotonic convergence under suitable conditions.
Notably, the additional condition here is readily verifiable, as it is on the initial condition and the heat source. 

\begin{theorem} \label{IBPpro}
{\bf (a)} The functions $w_m$, $u_m$ and $u_{m-1}$ satisfy the system:
\begin{align}\label{IBwmpro}
w_m (t,{\bf x},i) &= u_{m-1} (t,{\bf x},i) + \mathbb{E}_Q^{t,{\bf x},i} \left[ \mathbbm{1}(\tau_t^{(m)}\leq T\land\eta_t)\Theta_{t,\tau_t^{(m)}}w_0(\tau_t^{(m)},X_{\tau_t^{(m)}},R_{\tau_t^{(m)}}) \right],\\
u_m (t,{\bf x},i) &= u_{m-1} (t,{\bf x},i) + \mathbb{E}_Q^{t,{\bf x},i} \left[ \mathbbm{1}(\tau_t^{(m)}\leq T\land\eta_t)\Theta_{t,\tau_t^{(m)}}u_0(\tau_t^{(m)},X_{\tau_t^{(m)}},R_{\tau_t^{(m)}}) \right].\nonumber %label{umpro}
\end{align}

\noindent {\bf (b)}
It holds that for $(t,{\bf x},i) \in [0,T)\times D \times \mathcal{M}$, $|u_m(t,{\bf x},i)-v(t,{\bf x},i)|=\mathcal{O}((c^{m+1}(T-t)^{m+1}/(m+1)!)^p)$ as $m\to +\infty$ for any $p\in (0,1)$.
If each of the sequences $\{g(\cdot,i)\}_{i\in \mathcal{M}}$ and $\{\phi(\cdot,\cdot,i)\}_{i\in\mathcal{M}}$ is either non-negative or non-positive, then it holds that $u_m(\cdot,\cdot,i)\to  v(\cdot,\cdot,i)$ locally uniformly on $[0,T]\times D$ for $i\in\mathcal{M}$.

\noindent {\bf (c)}
If, moreover, $\{g(\cdot,i)\}_{i\in \mathcal{M}}$ is non-negative (resp. non-positive) and $\{\phi(\cdot,\cdot,i)\}_{i\in\mathcal{M}}$ is non-positive (resp. non-negative), then the locally uniform convergence is monotonic from below $u_m(\cdot,\cdot,i)\uparrow  v(\cdot,\cdot,i)$ (resp. from above $u_m(\cdot,\cdot,i)\downarrow  v(\cdot,\cdot,i)$) as $m\to +\infty$.
\end{theorem}

%\subsection{Hard bounding functions}

Moreover, as in Theorem \ref{IVPB}, we derive hard bounds based on
\[
M_r(t,{\bf x},i):=\sum_{j\in\mathcal{M}}q_{ij}({\bf x})\mathbb{E}_0^{t,{\bf x},j}\left[\int^{\eta_t\land T}_t \Theta_{t,s}^j ds\right],
\]
with a difference from \eqref{def of Mr} in the upper limit of the integral $\eta_t\land T$, instead of $T$.
With exactly the same definitions of the bound function $N_m$, the supremum and infimum $(N_m^U,N_m^L)$ and $(M_r^U,M_r^L)$, respectively, as \eqref{bound functions}, \eqref{def of NmU NmL} and \eqref{def of MrU MrL}, a slight modification of Theorem \ref{IVPB} yields the following result.
%\begin{equation*}
%N_m(t,{\bf x},i):=\left\{
%\begin{aligned}{}
%&\sum_{j\in\mathcal{M}}q_{ij}({\bf x})w_0(t,{\bf x},j), &\text{ if } m=0,\\
%&\sum_{j\in\mathcal{M}\setminus\{i\}}q_{ij}({\bf x})\left(w_m(t,{\bf x},j)-w_{m-1}(t,{\bf x},j)\right), &\text{ if } m\in \mathbb{N};
%\end{aligned}\right.
%\end{equation*}
%and for each bounded discounting rate $r(t,{\bf x},i)$, we define
%Let  be the essential supremum and infimum, respectively, of $N_m$ and $M_r$,
%\[
%N_m^U:=\esssup_{(t,{\bf x},i)\in[0,T]\times D \times\mathcal{M}}\left(N_m(t,{\bf x},i)\right)_+; \quad N_m^L:=\essinf_{(t,{\bf x},i)\in[0,T]\times D\times\mathcal{M}}\left(N_m(t,{\bf x},i)\right)_-;
%\]
%\[
%M_r^U:=\esssup_{(t,{\bf x},i)\in[0,T]\times D\times\mathcal{M}}\left(M_r(t,{\bf x},i)\right)_+; \quad M_r^L:=\essinf_{(t,{\bf x},i)\in[0,T]\times D\times\mathcal{M}}\left(M_r(t,{\bf x},i)\right)_-.
%\]
\begin{theorem}\label{IBB}
Suppose $|N_m^U|+|N_m^L|<+\infty$ for all $m\in\mathbb{N}_0$.

\noindent {\bf (a)}
If $M_r^U<1$, then it holds that
\[ %begin{equation}\label{boundswm0}
\frac{N_m^L}{1-M_r^L}\mathbb{E}_{0}^{t,{\bf x},i} \left[\int_{t}^{\eta_t\land T}\Theta_{t,s}^i ds\right]  
\leq v(t,{\bf x},i) - w_m(t,{\bf x},i)
\leq \frac{N_m^U}{1- M_r^U}\mathbb{E}_{0}^{t,{\bf x},i} \left[\int_{t}^{\eta_t\land T}\Theta_{t,s}^i ds\right],\quad (t,{\bf x},i)\in [0,T]\times D\times \mathcal{M}.
\] %end{equation}

\noindent {\bf (b)}
For each $(t,{\bf x},i)\in[0,T]\times D\times\mathcal{M}$, we have $|N_m(t,{\bf x},i)| \to  0$ as $m\to +\infty$.
\end{theorem}

\section{Numerical illustrations}\label{section numerical example}

In this section, we present numerical results to examine our theoretical findings, especially, the convergence results.
In order to achieve this primary goal without digressing into external numerical techniques at each step, we take a rather simple problem setting on purpose, where semi-analytical expressions are available for the target solution, the iterative approximations and the upper and lower bounding functions, so as to ensure a guaranteed accuracy of the numerical approximation throughout.
To this end, we consider the regime-switching geometric Brownian motion:
\begin{equation}\label{regime-switching geometric brownian motion}
dX_t = (r_{R_t}-\alpha_{R_t})X_tdt+\sigma_{R_t} X_t dW_t,
\end{equation}
where $\{W_t: t\geq 0\}$ is the standard Brownian motion in $\mathbb{R}$, and for $i\in\mathcal{M}$, $r_{i}$, $\alpha_i\in\mathbb{R}$ and $\sigma_i\in(0,\infty)$.
Clearly, this problem setting falls in the underlying stochastic differential equation (\ref{XSDE}) with $b(t,x,i)=(r_i-\alpha_i)x$ and $\sigma(t,x,i)=\sigma_i x$.
We let the regime process $\{R_t: t\geq 0\}$ be homogeneous.
%, that is, all entries of the generator matrix are independent of the state variable $x$.

Consider the target solution $v(t,x,i)= \mathbb{E}_Q^{t,x,i}[e^{-\int_t^T r_{R_s} ds}g(X_T)]$ on $[0,T]\times (0,+\infty)\times \mathcal{M}$, with a regime-independent non-negative initial condition $g$ and the homogeneous heat source $\phi\equiv 0$.
%We further assume that the initial condition $g(\cdot,i)=g(\cdot)$ is non-negative and independent of the second argument $i$.
On the basis of \eqref{wmpro} and \eqref{umpro}, the approximate functions $\{u_m\}_{m\in\mathbb{N}_0}$ and $\{w_m\}_{m\in\mathbb{N}_0}$ can be written recursively, with the aid of the function
\[
V(x,r,\sigma,t,\alpha) :=\int_{\mathbb{R}}e^{-r t} g\left( x \exp\left[(r-\alpha-\sigma^2/2)t +\sigma\sqrt{t}z\right]\right)\frac{1}{\sqrt{2\pi}}e^{-\frac{z^2}{2}}dz,
\]
which can be obtained by numerical approximation with any arbitrary order of accuracy.

\noindent {\bf (a)} We start with $m=0$, that is,
\[ 
w_0(t,x,i) = V(x,r_i,\sigma_i,T-t,\alpha_i), \quad u_0(t,x,i) = e^{q_{ii}(T-t)} V(x,r_i,\sigma_i,T-t,\alpha_i).
\]

\noindent {\bf (b)}
With $m=1$, we have
\begin{align*}
w_1(t,x,i) &= u_0(t,x,i) + \sum_{j\in\mathcal{M}\setminus\{i\}} q_{ij} \int_t^T e^{q_{ii}(s-t)} V(x, r(s),\sigma(s),T-t, \alpha(s)) ds,\\
u_1(t,x,i) &= u_0(t,x,i) + \sum_{j\in\mathcal{M}\setminus\{i\}} q_{ij} \int_t^T e^{q_{ii}(s-t)+q_{jj}(T-s)} V(x, r(s),\sigma(s),T-t, \alpha(s)) ds,
\end{align*}
where the functions $r(s)$, $\sigma(s)$ and $\alpha(s)$ are defined by
\[
r(s):=\frac{r_i(s - t)+r_j(T-s)}{T-t},\quad \sigma^2(s) := \frac{\sigma_i^2(s - t)+\sigma_j^2(T-s)}{T-t},\quad \alpha(s):=\frac{\alpha_i(s - t)+\alpha_j(T-s)}{T-t}.
\]

\noindent {\bf (c)}
Next, with $m=2$, we have
\begin{align*}
w_2(t,x,i) &= u_1(t,x,i) + \sum_{j\in\mathcal{M}\setminus\{i\}} \sum_{k\in\mathcal{M}\setminus\{j\}}q_{ij} q_{jk} 
\int_t^T \int_{s_1}^T e^{q_{ii}(s_1-t)+q_{jj}(s_2-s_1)} V(x, r(s_1,s_2),\sigma(s_1,s_2),T-t, \alpha(s_1,s_2))ds_2 ds_1,\\
u_2(t,x,i) &= u_1(t,x,i) + \sum_{j\in\mathcal{M}\setminus\{i\}} \sum_{k\in\mathcal{M}\setminus\{j\}}q_{ij} q_{jk} 
\int_t^T \int_{s_1}^T e^{q_{ii}(s_1-t)+q_{jj}(s_2-s_1)+q_{kk}(T-s_2)} V(x, r(s_1,s_2),\sigma(s_1,s_2),T-t, \alpha(s_1,s_2)) ds_2ds_1,
\end{align*}
where the functions $r(s_1,s_2)$, $\sigma(s_1,s_2)$ and $\alpha(s_1,s_2)$ are defined by
\begin{gather*}
r(s_1,s_2):=\frac{r_i(s_1 - t)+r_j(s_2-s_1)+r_k(T-s_2)}{T-t}, \quad\sigma^2(s_1,s_2) := \frac{\sigma_i^2(s_1 - t)+\sigma_j^2(s_2-s_1)+\sigma_k^2(T-s_2)}{T-t},\\ \alpha(s_1,s_2):=\frac{\alpha_i(s_1 - t)+\alpha_j(s_2-s_1)+\alpha_k(T-s_2)}{T-t}.
\end{gather*}
For $m=3$ and beyond, the functions $w_m$ and $u_m$ can be represented in terms of the previous approximations in a similar manner.
For more detail of the derivation, we refer the reader to the Appendix.

%Alternatively, the Fubini theorem allows us to reduce the expression for $\{w_m\}_{m\in\mathbb{N}}$ in (\ref{IVPwm}) and $\{u_m\}_{m\in\mathbb{N}}$ in (\ref{umpro0}) to:
%\begin{align*}
%w_m(t,x,i)
%&=e^{q_{ii} (T-t)}w_0(t,x,i)+ \sum_{j\in\mathcal{M}\setminus\{i\}}q_{ij}  \int_t^T e^{(q_{ii}-r_i)(s-t)} \mathbb{E}_0^{t,x,i}\left[ w_{m-1}(s,X_s,j) \right] ds, \\
%u_m(t,x,i)
%&=e^{q_{ii} (T-t)}w_0(t,x,i)+ \sum_{j\in\mathcal{M}\setminus\{i\}}q_{ij}  \int_t^T e^{(q_{ii}-r_i)(s-t)} \mathbb{E}_0^{t,x,i}\left[ u_{m-1}(s,X_s,j) \right] ds.
%\end{align*}
%Under $\mathbb{P}^{t,{\bf x},i}_0$, we have the equality in distribution
%$
%X_s \stackrel{\mathcal{L}}{=} x \exp\left((r_i-\alpha_i-\sigma_i^2/2)(s-t)+(\sigma_i\sqrt{s-t})Z\right)
%$,
%where $Z$ is a standard normal random variable, thus the expectation in the above equation can be computed as an integral with the probability density function of the log-normal distribution.

Next, in light of Section \ref{section hard bounding functions IVP}, we derive upper and lower hard bounding functions.
Note first $\mathbb{E}_{0}^{t,x,i}[\int_{t}^{T} e^{- r_{R_s}(s-t) } ds] = (1-e^{-r_i (T-t)})/r_i$, due to $R_s\equiv i$ for all $s\in [t,T]$, $\mathbb{P}_{0}^{t,x,i}$-$a.s.$
Hence, one can compute the bounds for the target solution at $m$-th iteration on the basis of (\ref{LUm}):
\begin{equation}\label{LUmsimple2}
U_m(t,x,i;\,w)=w_m(t,x,i)+\frac{N_m^U(w)}{1-M_r^U}\left( \frac{1-e^{-r_i (T-t)}}{r_i}\right),\quad 
L_m(t,x,i;\,w)=w_m(t,x,i)+\frac{N_m^L(w)}{1-M_r^L}\left( \frac{1-e^{-r_i (T-t)}}{r_i}\right),
\end{equation}
where $M_r^U$ and $M_r^L$ are respectively the essential infimum and supremum of the function $M_r(t,x,i) = \sum_{j\in\mathcal{M}}q_{ij} (1-e^{-r_j (T-t)})/r_j$, which is independent of the state variable $x$, and $N_m^L(w)$ and $N_m^U(w)$ are respectively the essential infimum and supremum of \eqref{bound functions}.
We obtain $\{N_m^U(w)\}_{m\in\mathbb{N}_0}$, $\{N_m^L(w)\}_{m\in\mathbb{N}_0}$, $M_r^U$ and $M_r^L$ by numerical approximation, as it seems difficult to obtain those here in an analytic manner.

For a clearer and thorough demonstration, we begin with a two-regime setting with $g(\cdot,i)=(\cdot-K)_+$ for some $K>0$, so that a direct comparison with relevant existing methods \cite{louis, optionP} is possible.
Now, we start with the initial approximation $w_0(t,x,i)=C(x,r_i,\sigma_i,T-t,\alpha_i,K)$, where 
\begin{equation*}
C(x,r,\sigma,t,\alpha,K) := x e^{-\alpha t}\Phi\left(\frac{\ln(x/K)+(r-\alpha+\sigma^2/2)t}{\sigma \sqrt{t}}\right)
-Ke^{-rt} \Phi\left(\frac{\ln(x/K)+(r-\alpha-\sigma^2/2)t}{\sigma \sqrt{t}}\right),
\end{equation*}
where we denote by $\Phi$ the standard normal cumulative distribution function.
We set the parameters to $q_{12}=q_{21}=1.0$, $r_1=r_2=0.05$, $\sigma_1=0.15$, $\sigma_2=0.25$, $\alpha_1=\alpha_2=0$, $T=1.0$ and $K=1.0$.
It is worth mentioning that due to the homogeneous setting $r_1=r_2$, the upper and lower bounding functions \eqref{LUmsimple2} further simplifies according to \eqref{LUmsimple}.
Note also that we are here expecting a uniform and monotonic convergence of the sequence $\{u_m(\cdot,\cdot,i)\}_{m\in\mathcal{M}}$ since the initial conditions $\{g(\cdot,i)\}_{i\in\mathcal{M}}$ are non-negative and the heat sources $\{\phi(\cdot,i)\}_{i\in\mathcal{M}}$ are flat zero (Theorem \ref{IVPpro} {\bf (b)}).

\begin{figure}[ht] 
  \centering
  \begin{subfigure}{0.38\linewidth}
    \includegraphics[width=\linewidth]{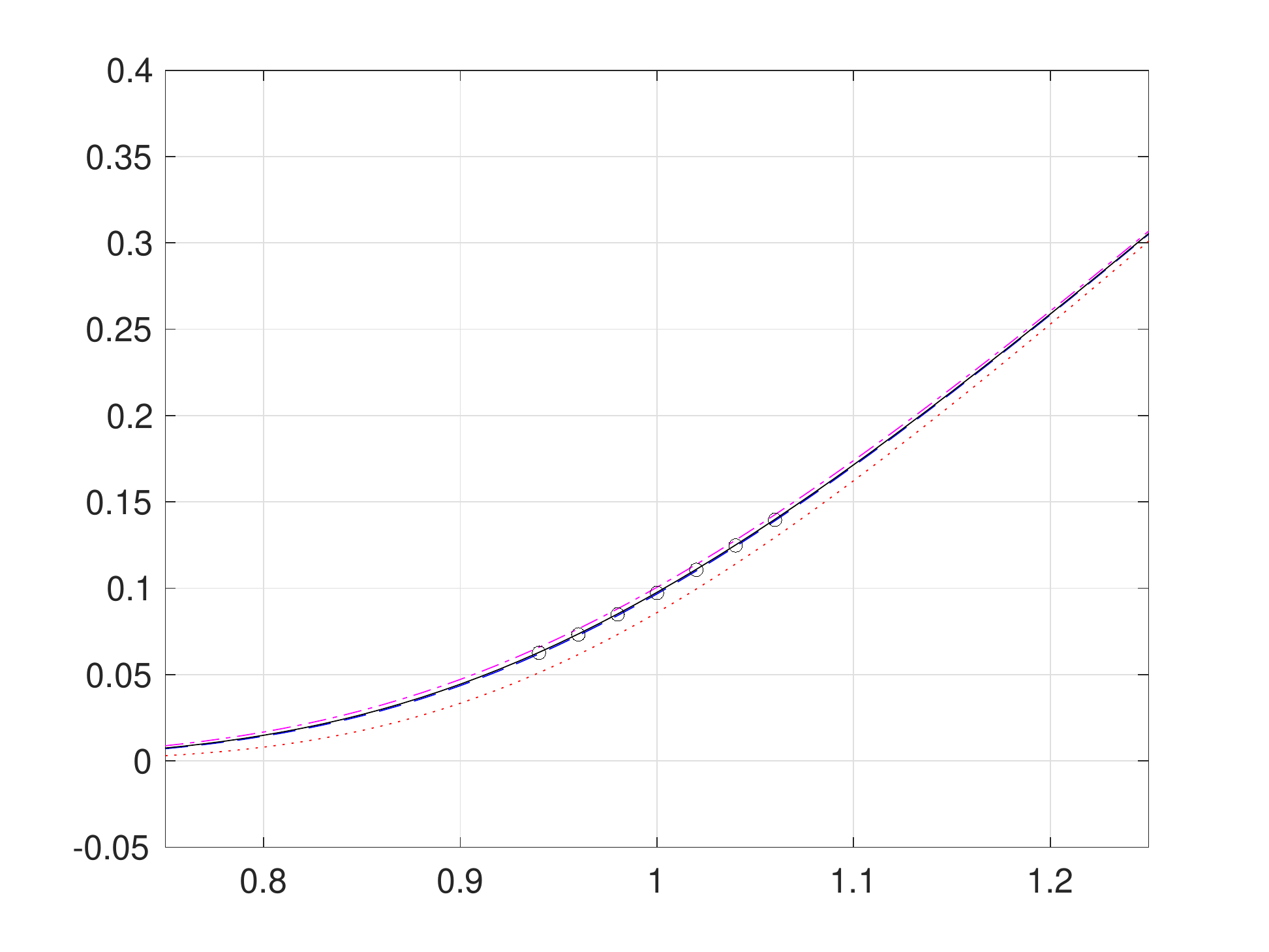}
    \caption{Iterative approximations $w_m(0,\cdot,1)$}
  \end{subfigure}
  \begin{subfigure}{0.38\linewidth}
    \includegraphics[width=\linewidth]{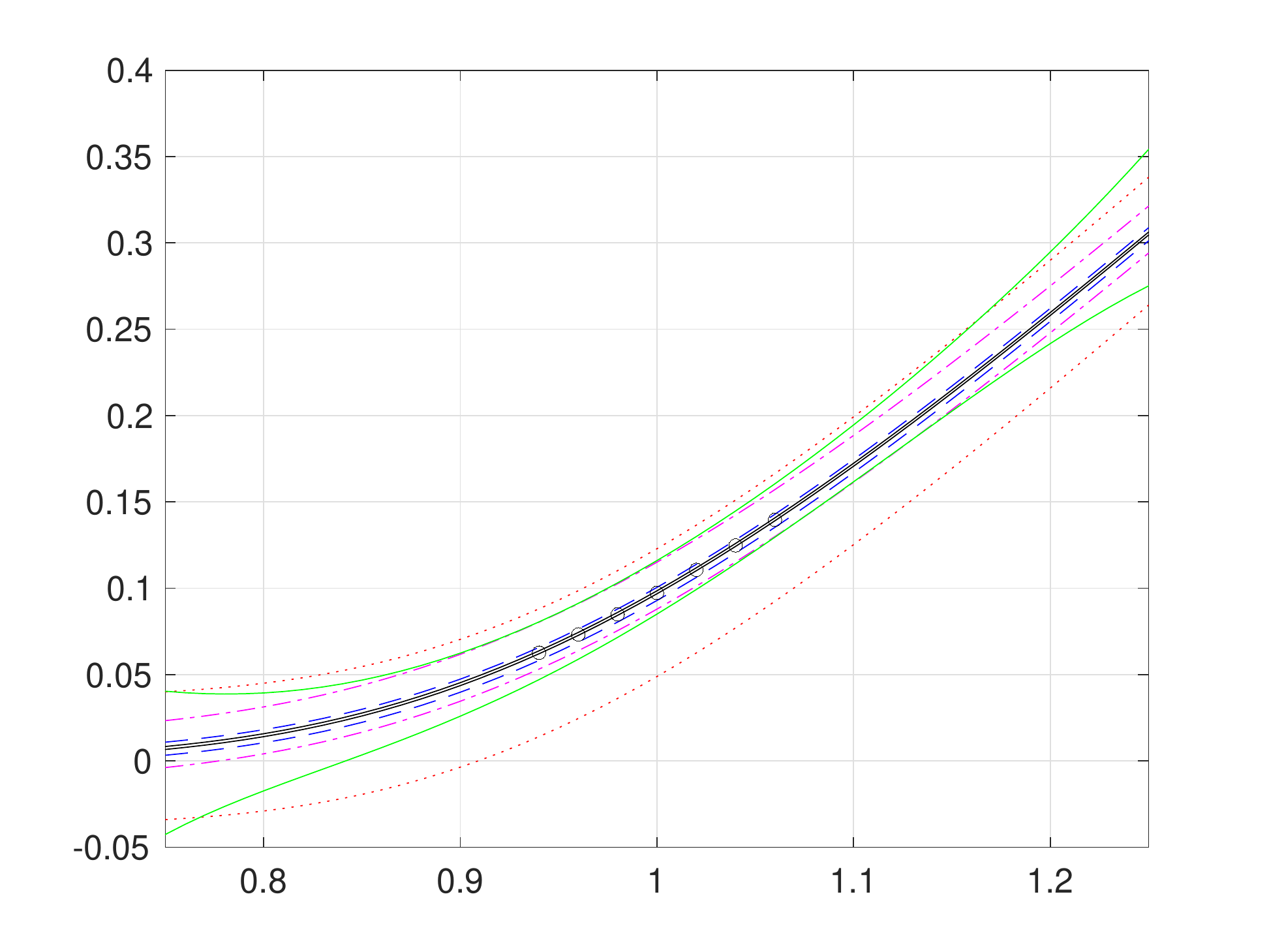}
    \caption{Iterative bounds $L_m(0,\cdot,1;\,w)$ and $U_m(0,\cdot,1;\,w)$}
  \end{subfigure}
  
\begin{subfigure}{0.38\linewidth}
    \includegraphics[width=\linewidth]{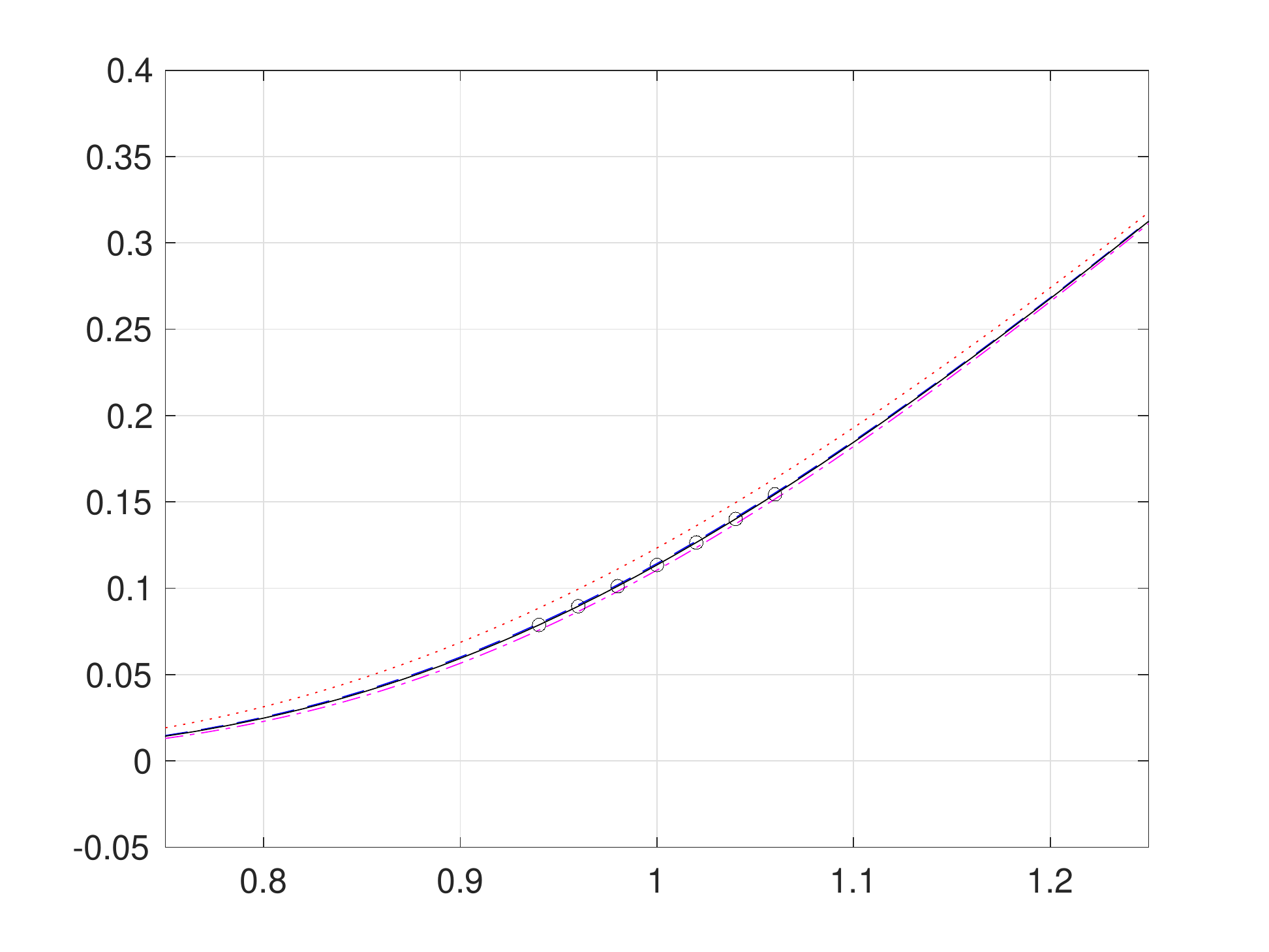}
    \caption{Iterative approximations $w_m(0,\cdot,2)$}
\end{subfigure}
  \begin{subfigure}{0.38\linewidth}
    \includegraphics[width=\linewidth]{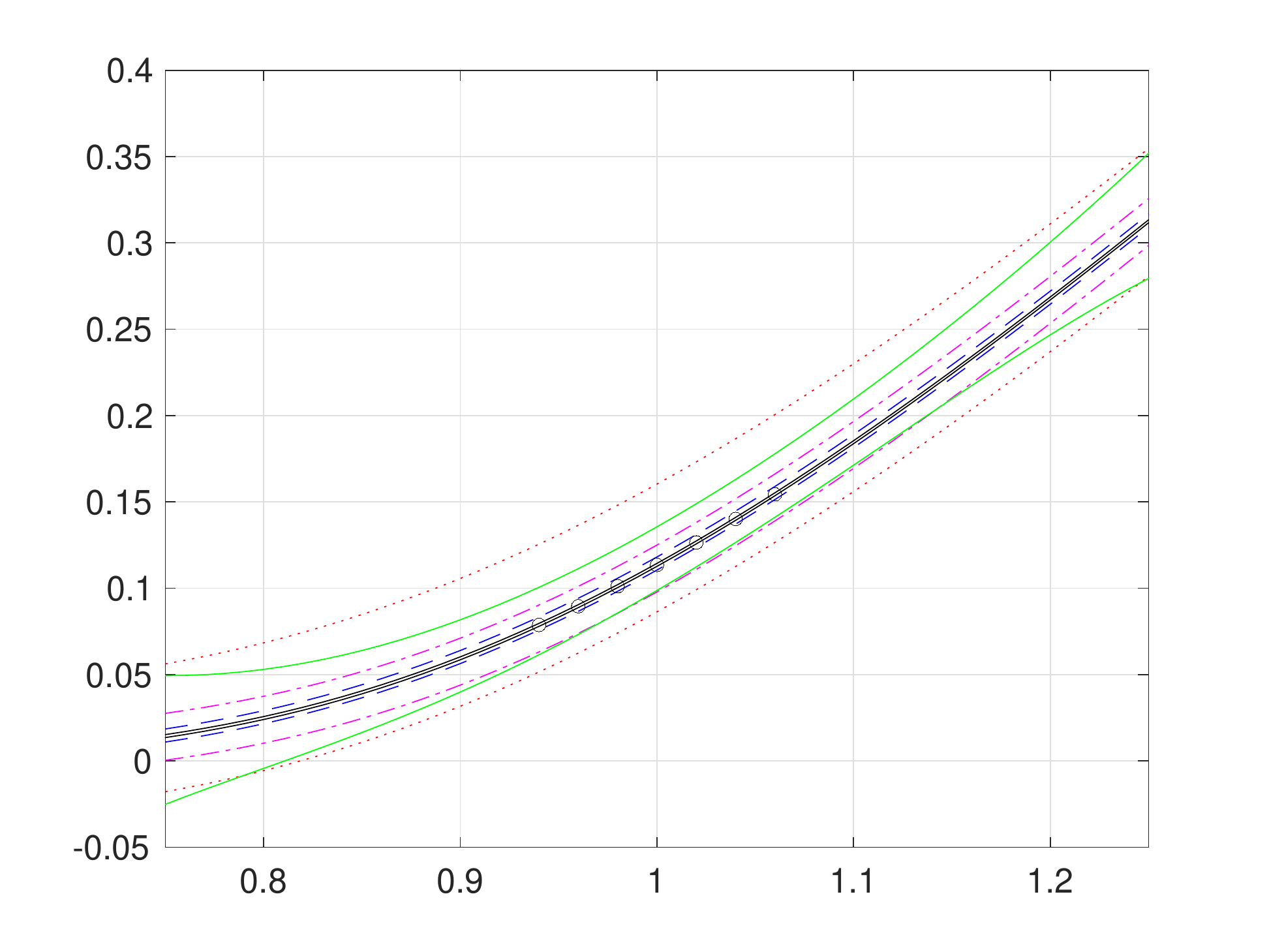}
    \caption{Iterative bounds $L_m(0,\cdot,2;\,w)$ and $U_m(0,\cdot,2;\,w)$}
  \end{subfigure}
  \caption{Numerical results for the two-regime setting.
  Each figure contains $m=0$ (red dot), $m=1$ (pink dash-dot), $m=2$ (blue dash) and $m=3$ (black solid).
The circles provide pointwise approximations as reported in \cite{optionP}, while the green solid lines  in (b) and (d) provide the polynomial bounds \cite{louis} of degree 10 (green solid).}
  \label{fig01}
\end{figure}

First, we present in Figure \ref{fig01} iterative weak approximations $\{w_m(0,\cdot,i)\}_{m\in \{0,1,2,3\}}$, along with sequences of hard bounding functions $\{L_m(0,\cdot,i;\,w)\}_{m\in \{0,1,2,3\}}$ and $\{U_m(0,\cdot,i;\,w)\}_{m\in\{0,1,2,3\}}$ for the first three iterations $m\in\{0,1,2,3\}$.
For a direct comparison with existing numerical results in the literature, we also provide the pointwise approximations reported in \cite{optionP}, as well as the polynomial upper and lower bounds of \cite{louis}.
Our hard bounding functions converge towards each other quite fast, largely because it is highly likely to observe switching at most three times within the interval $[0,T]$, indeed with a $98.1\% (=\mathbb{P}_Q^{0,{\bf x},i} (\tau_0^{(4)} >T) = (1+1+1/2+1/6) e^{-1})$ chance. 

\begin{figure}[ht] 
  \centering
  \begin{subfigure}{0.33\linewidth}
    \includegraphics[width=\linewidth]{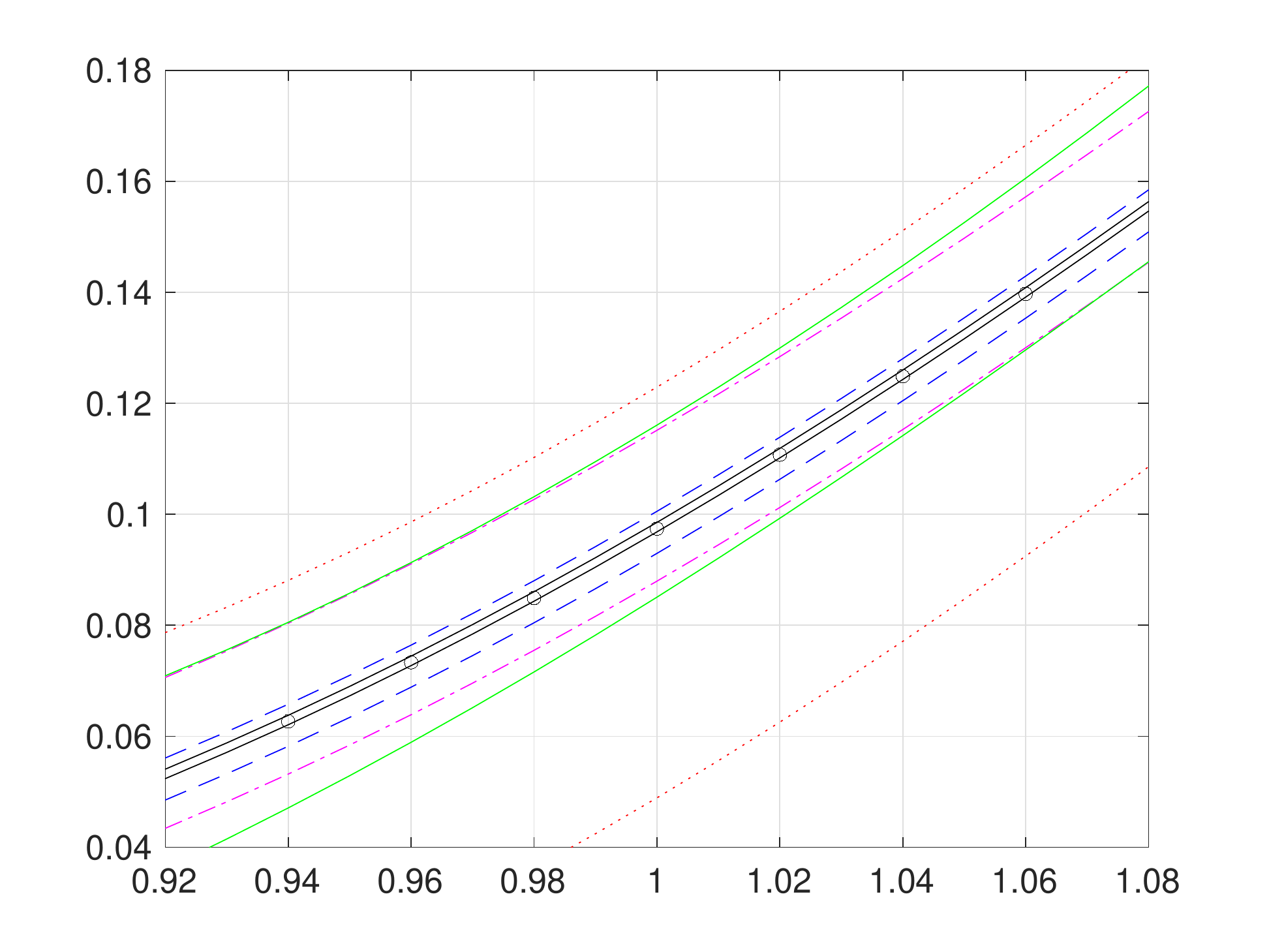}
    \caption{Zoom-in of Figure \ref{fig01} (b)}
\end{subfigure}
  \begin{subfigure}{0.33\linewidth}
    \includegraphics[width=\linewidth]{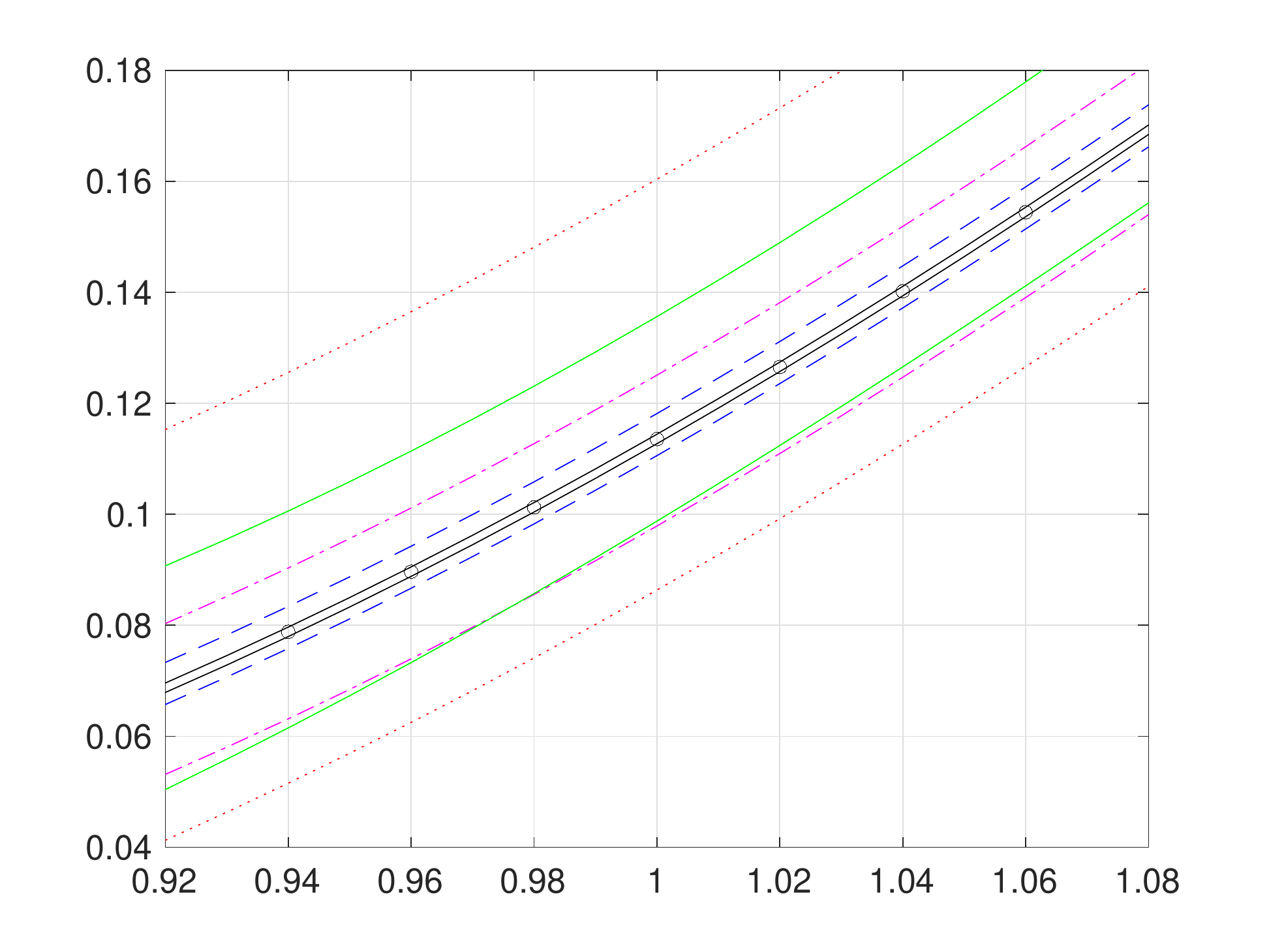}
    \caption{Zoom-in of Figure \ref{fig01} (d)}
    \end{subfigure}
  \begin{subfigure}{0.33\linewidth}
    \includegraphics[width=\linewidth]{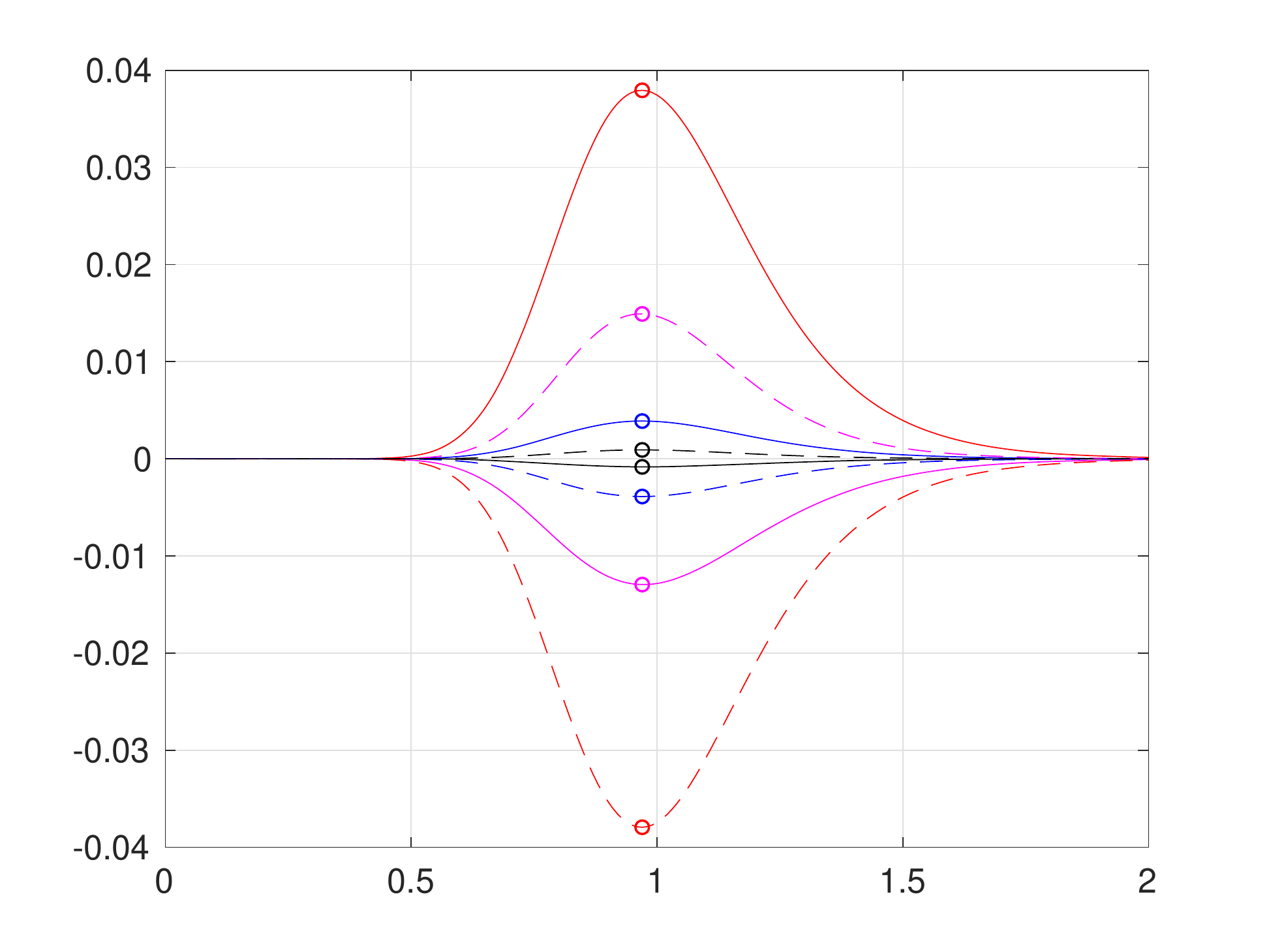}
    \caption{Bound functions $N_m(0,\cdot,i;w)$ of \eqref{bound functions} for $m\in\{0,1,2,3\}$ and $i\in \{1,2\}$}
  \end{subfigure}
  \caption{Numerical results for the two-regime setting.
  Figures (a) and (b) zoom in Figure \ref{fig01} (b) and (d), respectively.
  Figure (a) provides the bound functions $N_m(0,\cdot,i)$ for $i=1$ (solid) and $i=2$ (dash), with $m=0$ (red), $m=1$ (pink), $m=2$ (blue) and $m=3$ (black).}
  \label{fig02}
\end{figure}

Figure \ref{fig01} is overly congested with lines and circles to fully illustrate strong agreement of our convergent hard bounding functions $\{L_m(0,\cdot,i;\,w)\}_{m\in\{0,1,2,3\}}$ and $\{U_m(0,\cdot,i;\,w)\}_{m\in\{0,1,2,3\}}$ with the existing pointwise approximations reported in \cite{optionP}. 
Hence, for a better presentation, we provide Figure \ref{fig02} (a) and (b), respectively, zoom-ins of Figure \ref{fig01} (b) and (d). 
In addition, Figure \ref{fig02} (c) presents the bounds function $N_m(0,x,i)$ defined in \eqref{bound functions}, from which the supremum and infimum $\{N_m^L(w)\}_{m\in\{0,1,2,3\}}$ and $\{N_m^U(w)\}_{m\in\{0,1,2,3\}}$ are obtained by numerical approximation, as follows: 
\begin{multline*}
(N_0^L(w),N_0^U(w), N_1^L(w), N_1^U(w), N_2^L(w), N_2^U(w), N_3^L(w), N_3^U(w))\\
 \approx (-0.0379,+0.0379, -0.0129,+0.0149, -0.0039,+0.0039, -0.00083,+0.00091).
\end{multline*}

\begin{figure}[ht]
  \centering
  \begin{subfigure}{0.35\linewidth}
    \includegraphics[width=\linewidth]{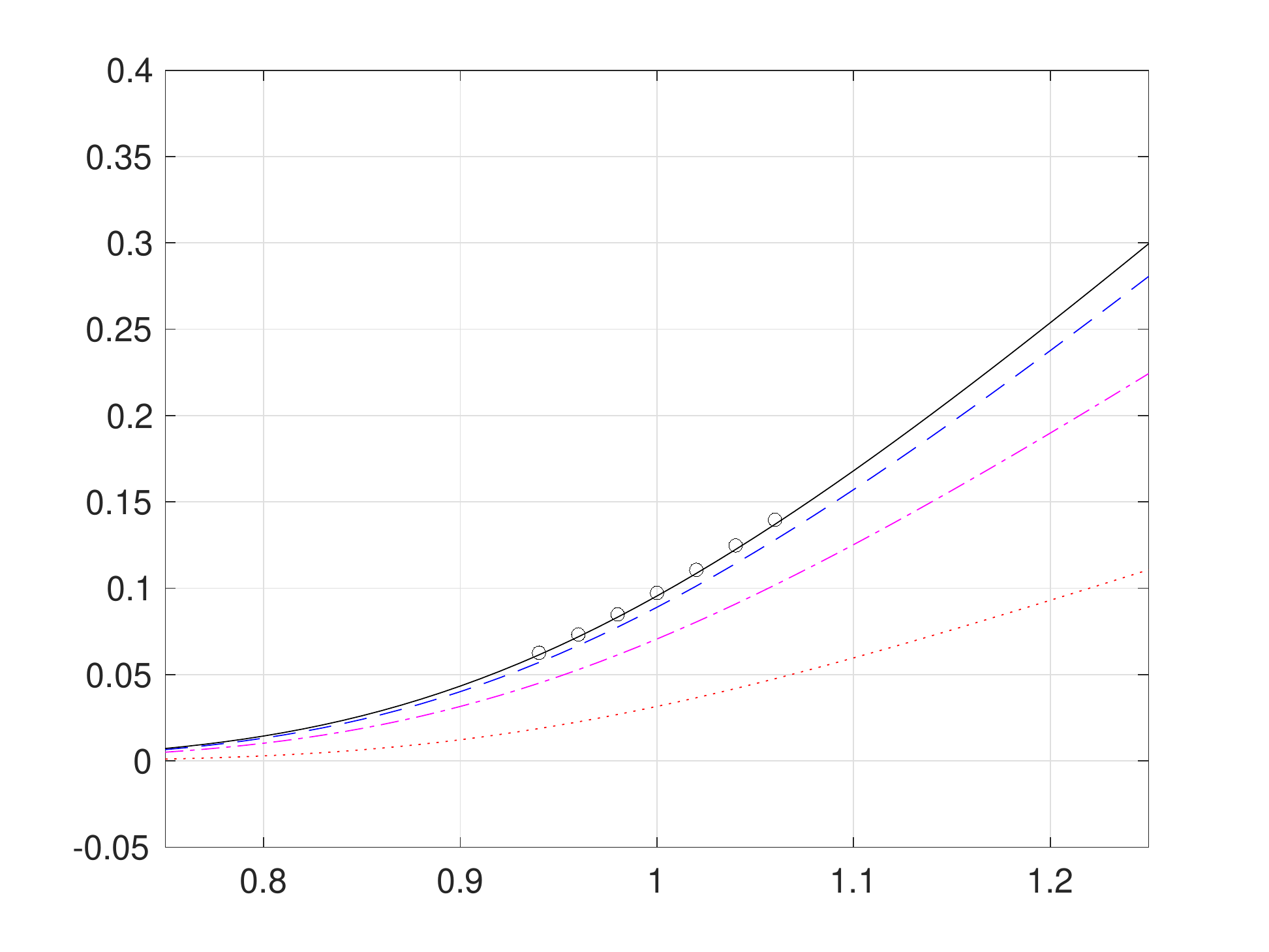}
     \caption{Regime 1: lower bounds $u_m(0,\cdot,1)$}
  \end{subfigure}
\begin{subfigure}{0.35\linewidth}
    \includegraphics[width=\linewidth]{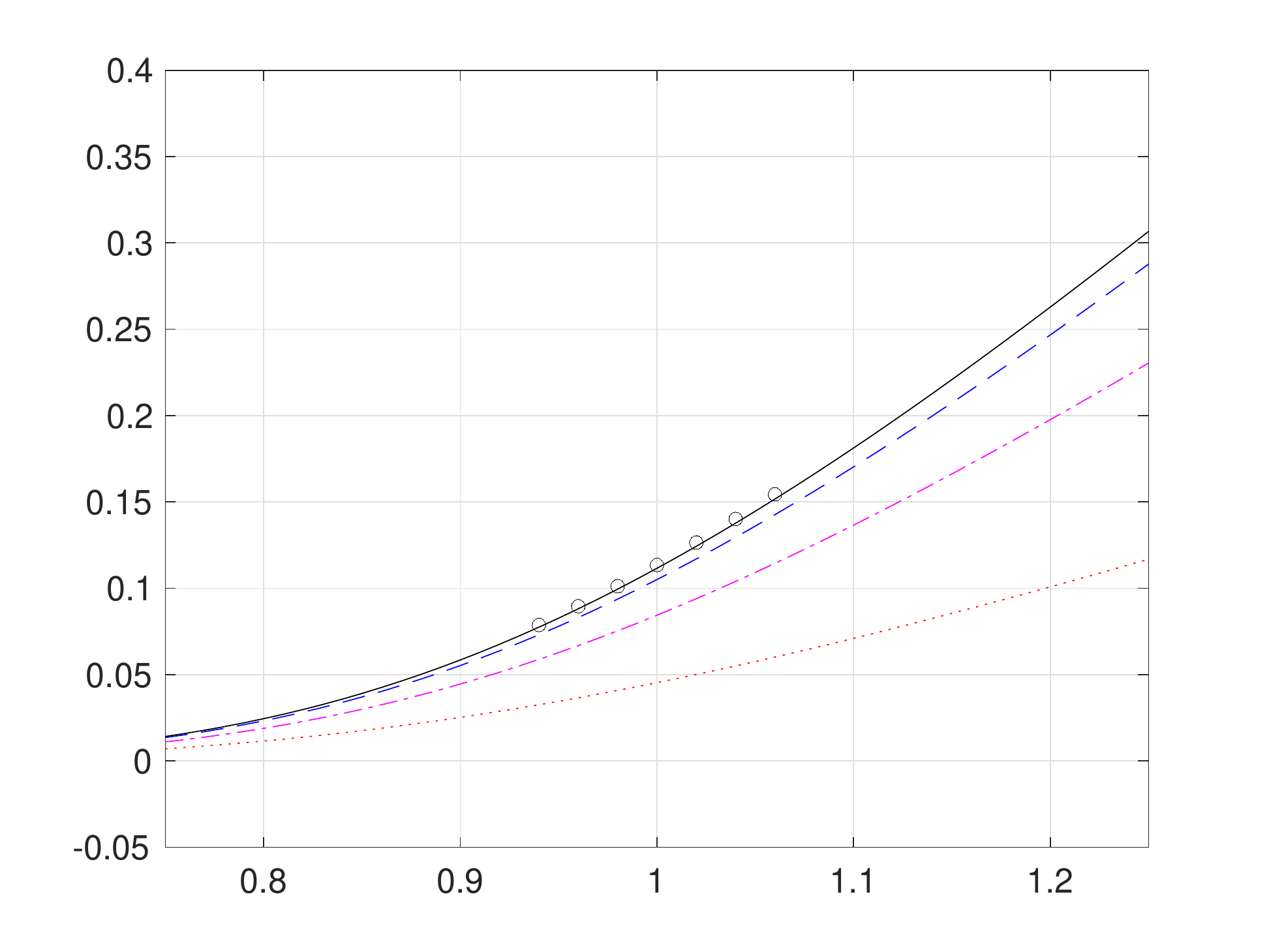}
    \caption{Regime 2: lower bounds $u_m(0,\cdot,2)$}
\end{subfigure}
 \caption{Numerical results for monotonic convergence of $\{u_m\}_{m\in\mathbb{N}}$ in the two-regime setting.
  Each figure contains $m=0$ (red dot), $m=1$ (pink dash-dot), $m=2$ (blue dash) and $m=3$ (black solid).
  The circles provide point-wise approximations as reported in \cite{optionP}.
%(a) Lower bounding functions $\{u_m(0,x,1)\}_{m\in\{0,1,2,3\}}$, along with the point-wise approximations (black circles) reported in \cite{optionP}; (b) is the one starting from regime 2, corresponding to (a).
}
\label{fig03}
\end{figure}

\begin{figure}[ht] 
  \centering
  \begin{subfigure}{0.33\linewidth}
    \includegraphics[width=\linewidth]{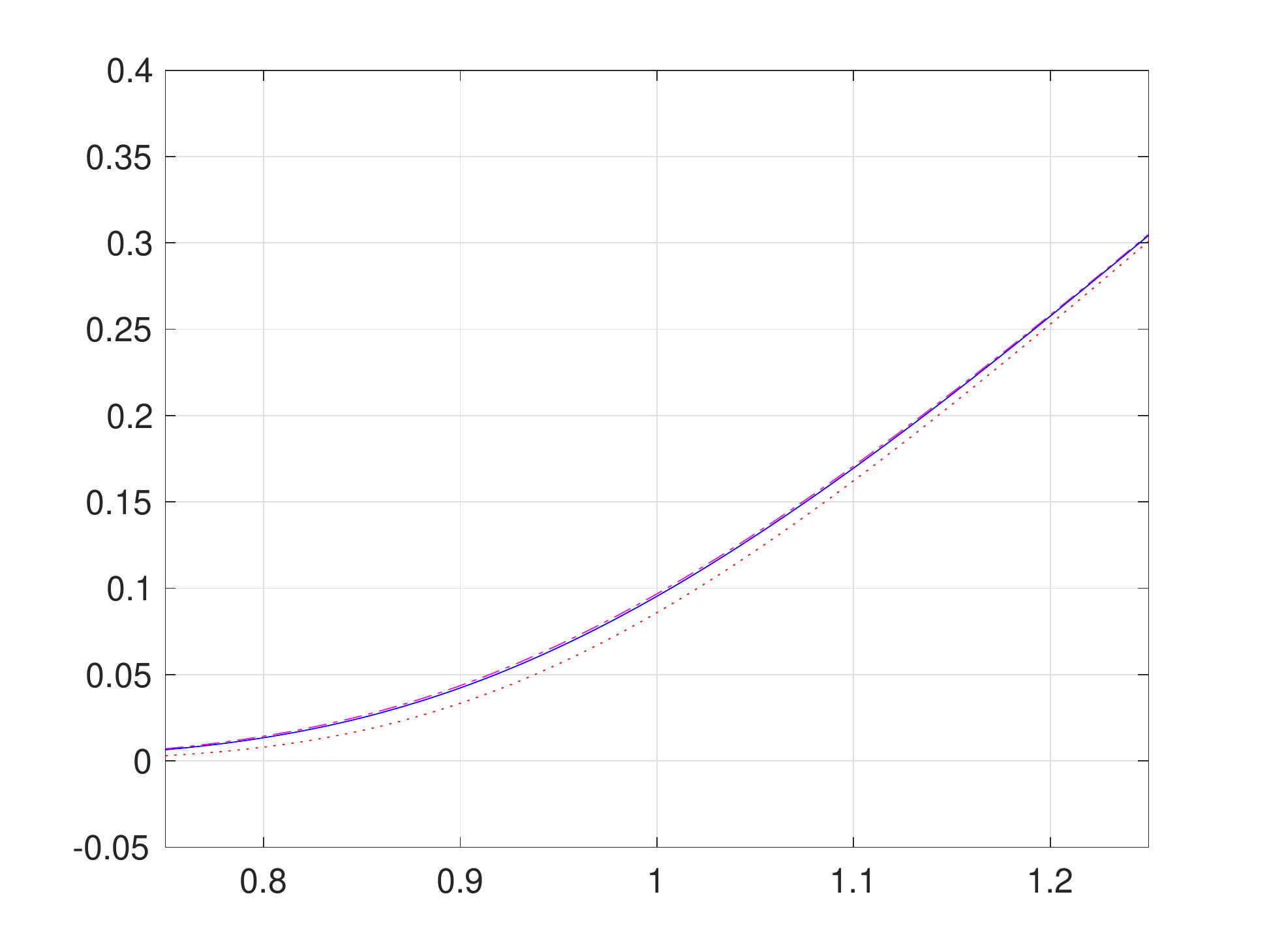}
    \caption{Iterative approximations $w_m(0,\cdot,1)$}
\end{subfigure}
  \begin{subfigure}{0.33\linewidth}
    \includegraphics[width=\linewidth]{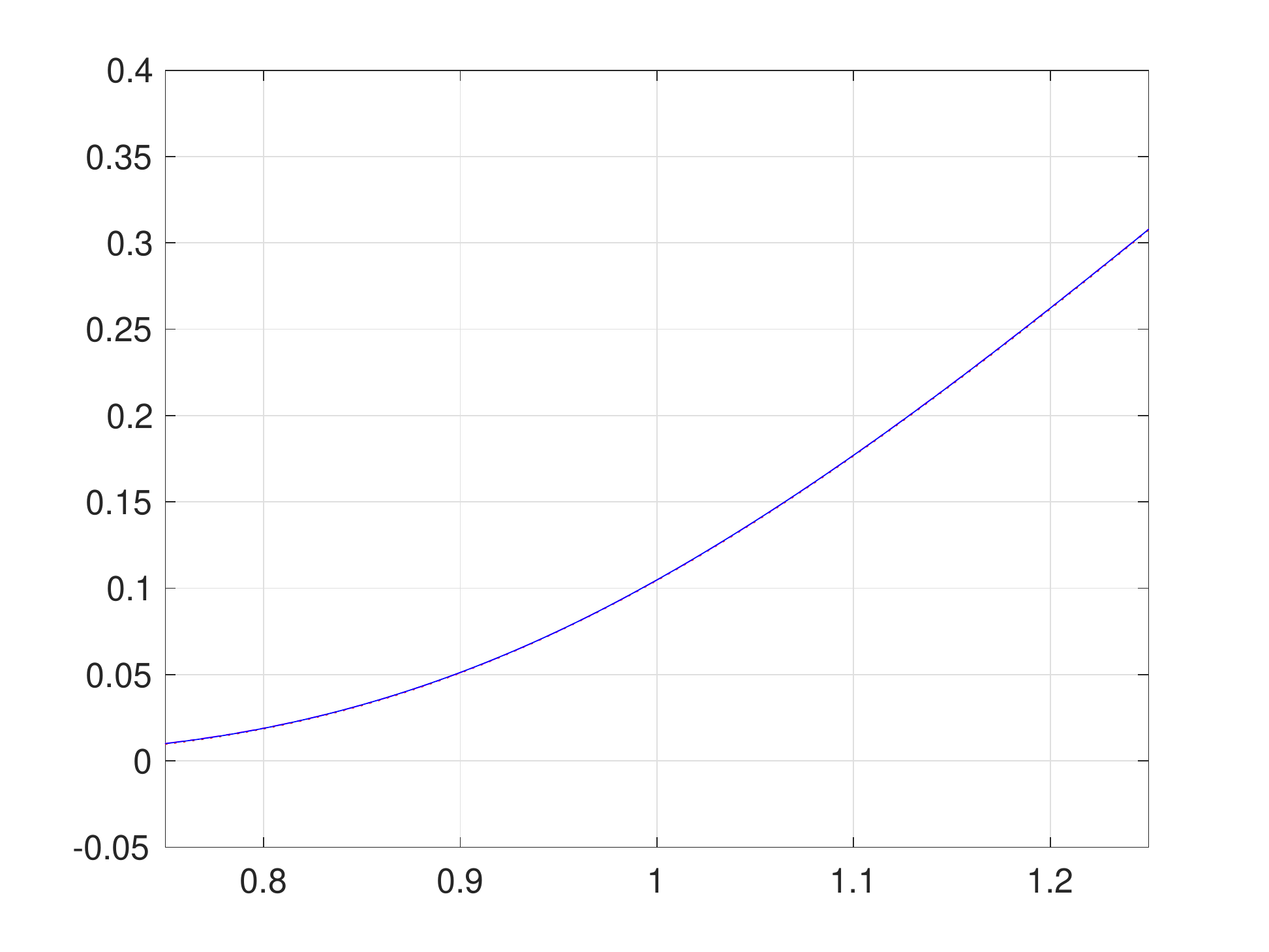}
    \caption{Iterative approximations $w_m(0,\cdot,2)$}
\end{subfigure}
  \begin{subfigure}{0.33\linewidth}
    \includegraphics[width=\linewidth]{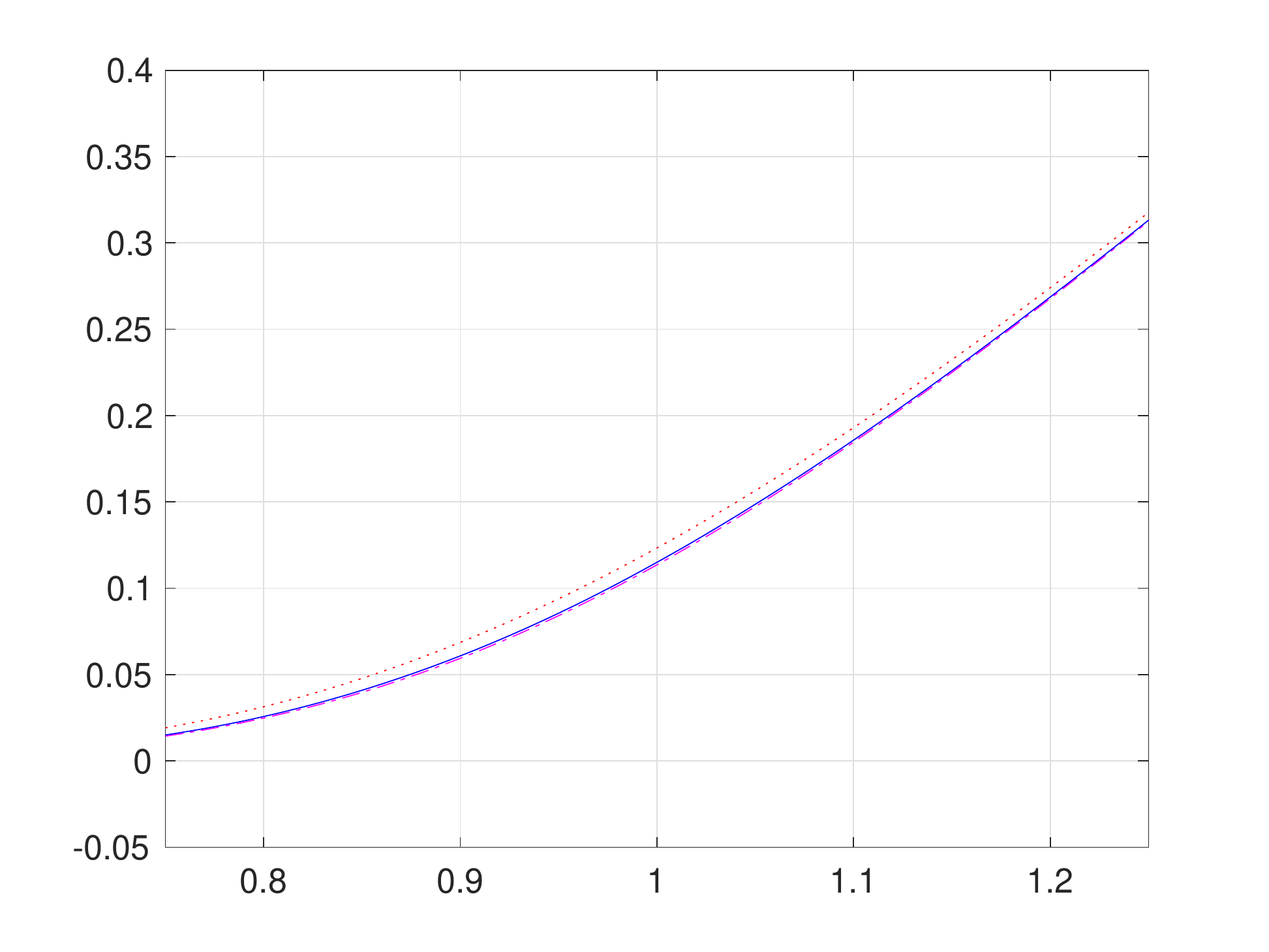}
    \caption{Iterative approximations $w_m(0,\cdot,3)$}
  \end{subfigure}
    \begin{subfigure}{0.33\linewidth}
    \includegraphics[width=\linewidth]{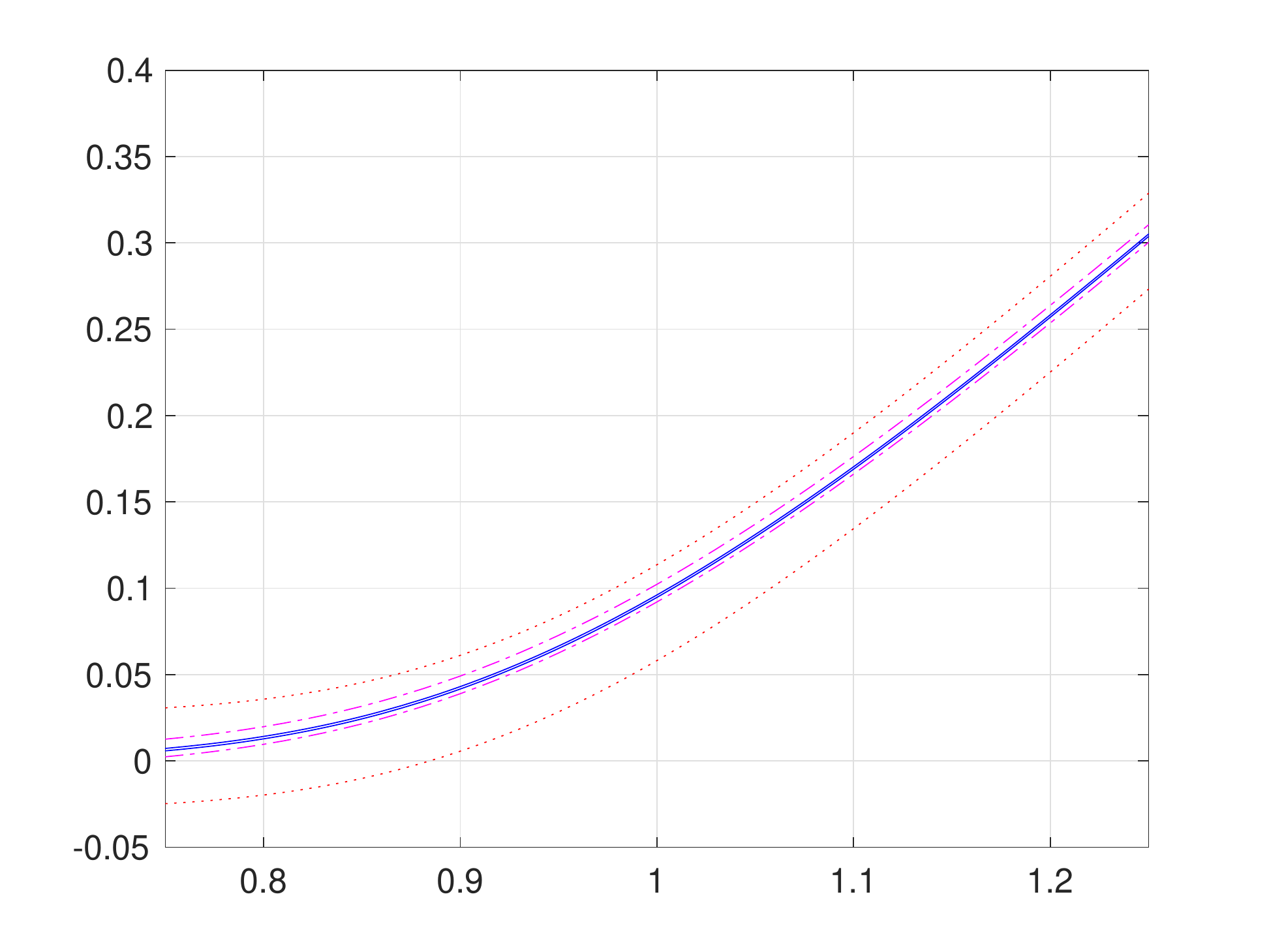}
    \caption{Iterative bounds $L_m(0,\cdot,1)$ and $U_m(0,\cdot,1)$}
  \end{subfigure}
    \begin{subfigure}{0.33\linewidth}
    \includegraphics[width=\linewidth]{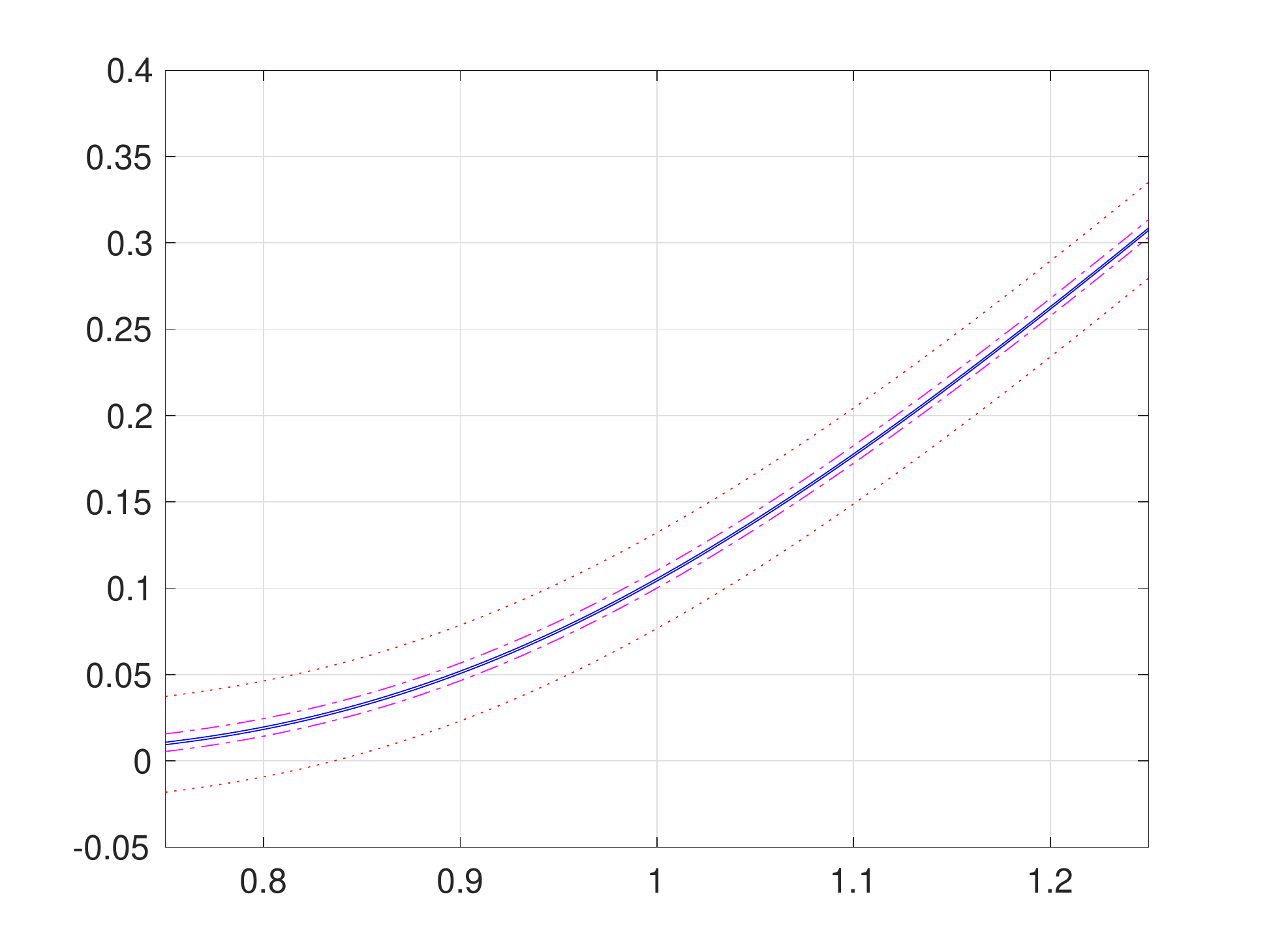}
    \caption{Iterative bounds $L_m(0,\cdot,2)$ and $U_m(0,\cdot,2)$}
  \end{subfigure}
      \begin{subfigure}{0.33\linewidth}
    \includegraphics[width=\linewidth]{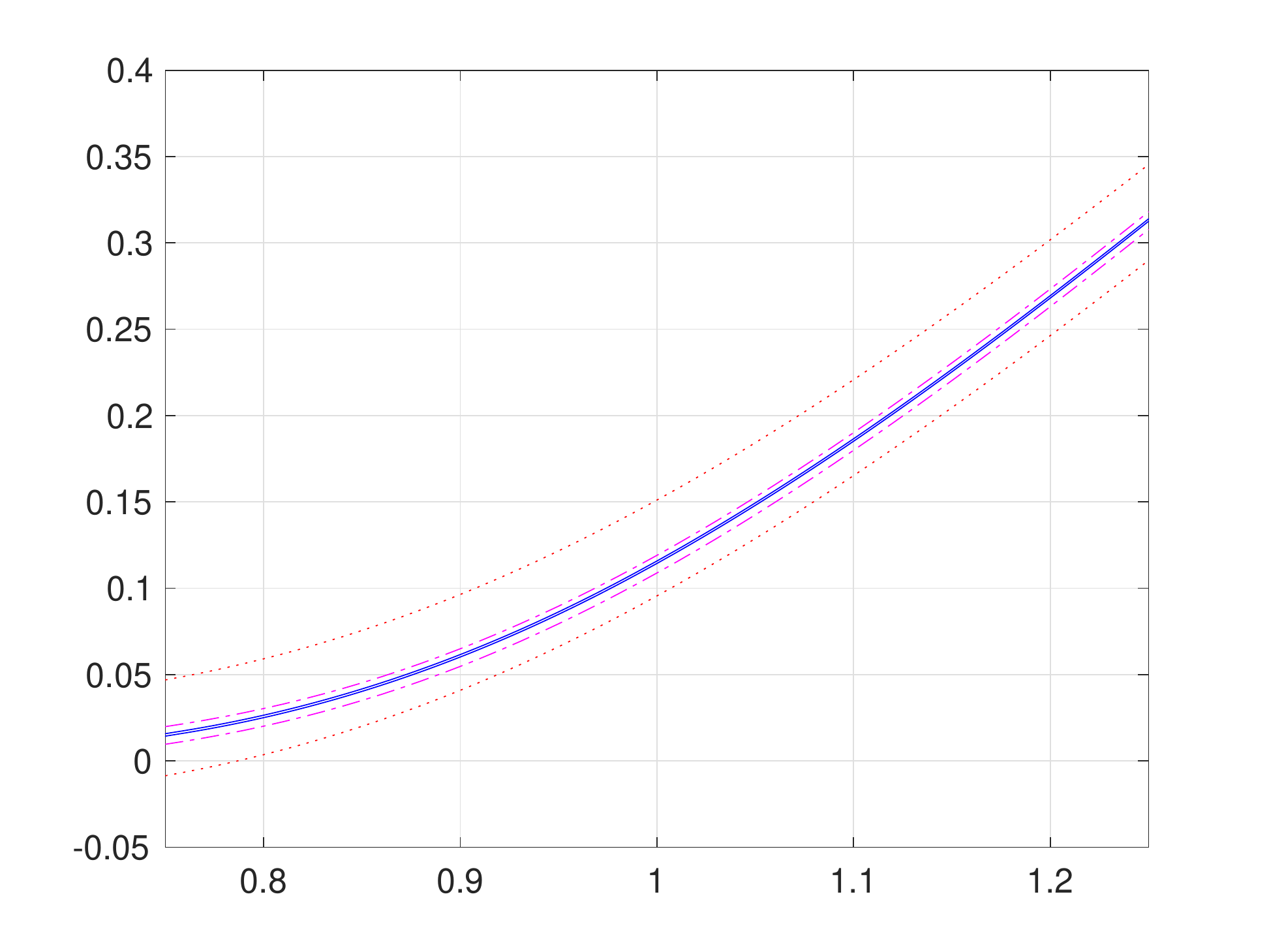}
    \caption{Iterative bounds $L_m(0,\cdot,3)$ and $U_m(0,\cdot,3)$}
  \end{subfigure}
 %       \begin{subfigure}{0.33\linewidth}
%    \includegraphics[width=\linewidth]{Nm12.eps}
 %   \caption{bound function $N_m(0,x,i)$}
 % \end{subfigure}
    \caption{Numerical results for the three-regime setting with $q_{ij}=0.5$ for all $i\ne j$, $r_1=r_2=r_3=0.05$, $(\sigma_1,\sigma_2,\sigma_3) = (0.15,0.2,0.25)$, $\alpha_1=\alpha_2=\alpha_3=0$, $T=1$ and $K=1.$
    Each figure contains $m=0$ (red dot), $m=1$ (pink dash-dot) and $m=2$ (blue solid).
%(a) Approximate functions $\{w_m(0,x,1)\}_{m\in\{0,1,2\}}$; (d) Corresponding bounding functions at $m$-th iteration $\{L_m(0,x,1),U_m(0,x,1)\}_{m\in\{0,1,2\}}$; Figures (b) and (e) draw numerical results when the dynamic starts from regime-2; Figures (c) and (f) draw numerical results when the dynamic starts from regime-3;  %Figure (g) Bound functions $N_m(0,x,i)$ for $i=1$ (solid), $i=2$ (dash) and $i=3$ (dash-dot), with $m=0$ (red), $m=1$ (pink) and $m=2$ (blue).
}
  \label{fig04}
\end{figure}

In this problem setting, the initial conditions $\{g(\cdot,i)\}_{i\in\mathcal{M}}$ are non-negative and the heat sources $\{\phi(\cdot,i)\}_{i\in\mathcal{M}}$ are flat zero.
Hence, by Theorem \ref{IVPpro} {\bf (b)} and {\bf (c)}, we are here expecting a uniform convergence of the sequence $\{u_m(\cdot,\cdot,i)\}_{m\in\mathbb{N}}$ to be also monotonic from below of the target solution $v$.
In Figure \ref{fig03}, we present the first three iterations, that is, $u_m(t,x,i)$ for $m\in\{0,1,2,3\}$.
As expected, the iterative approximations converge monotonically from below (as well as very quickly) towards the pointwise approximations as reported in \cite{optionP}.
%To avoid overloading the figures and for clearer presentation, we do not provide numerical results for Theorem \ref{IVPB} {\bf (c)}. 
We close this section with Figure \ref{fig04} for numerical results of the sequence $\{w_m\}_{m\in\mathbb{N}}$ and associated hard bounds $\{U_m(\cdot,\cdot,\cdot;\,w)\}_{m\in\mathbb{N}_0}$ and $\{L_m(\cdot,\cdot,\cdot;\,w)\}_{m\in\mathbb{N}_0}$ for a three-regime case.

%\textcolor{red}{Moreover, the sequence $\{u_m\}$ is monotonically increasing in $x$ in this example, thus it can also be interpreted as an alternative lower bound for $v(0,{\bf x},i)$.}

%{\small
%}

\section{Concluding remarks}\label{section concluding remarks}

In this paper, we have established a novel convergent iteration framework for a weak approximation of general switching diffusion.
By restricting the maximum number of switching, we succeeded to untangle a challenging system of weakly coupled partial differential equations to a collection of independent partial differential equations, for which a variety of accurate and efficient numerical methods are available.
Using the sequence of resulting approximate solutions, we constructed upper and lower bounding functions for the unknown solutions.
Convergences of the iterative approximate solutions as well as the associated upper and lower bounding functions were rigorously derived and demonstrated through numerical results.

We believe that the developed framework has the potential to be a viable tool for the weak approximation and hard bounding functions for a general class of stochastic differential equations with regime switching.
In the present work, in order to fully carry out our theoretical developments and demonstrate all aspects of the developed framework through numerical results, we restricted numerical illustrations to a minimum and did not go into a study of its range of applicability, which is significantly large and should thus be demonstrated in the future.
Other future research topics include a investigation of possible issues and cumulative error caused by an iterative use of external numerical techniques, for instance, finite difference methods to approximate $w_m$ based on the previous approximation for $w_{m-1}$ in the system \eqref{IVPwmpde}, in light of various factors, such as the problem dimension and the input data.

\small
\bibliographystyle{abbrv}
\bibliography{regimebib}

\appendix

\section{Technical notes}\label{appendix technical notes}
For the sake of completeness, we provide some details of derivation for the computation conducted in the numerical experiments of Sections \ref{section numerical example}.
%We begin with the derivation for Sections \ref{section numerical example}.
Since no switching is allowed under $\mathbb{P}_0^{t,x,i}$, it holds that $\{X_s:s \geq t\}$ satisfies the stochastic differential equation $dX_s=(r_i-\alpha_i)X_s ds+\sigma_i X_s dW_s$ without switching, that is, under the probability measure $\mathbb{P}_0^{t,x,i}$,
\[
X_T  \stackrel{\mathcal{L}}{=} x \exp\left[(r_i-\alpha_i-\sigma_i^2/2)(T-t)+\sigma_i\sqrt{T-t} Z\right],
\]
where we reserve $Z$ for a general standard normal random variable, independent of the $\sigma$-field $\mathcal{F}$.
Hence, we have
\[
w_0(t,x,i) = \mathbb{E}_0^{t,x,i}\left[e^{-r_i(T-t)}g(X_T)\right] = \mathbb{E}\left[e^{-r_i(T-t)}g\left(x \exp\left[(r_i-\alpha_i-\sigma_i^2/2)(T-t)+\sigma_i\sqrt{T-t}Z\right]\right)\right]
=V(x,r_i,\sigma_i,T-t,\alpha_i).
\]
Since the transition rate $-q_{ii}$ is independent of the state variable $x$ in this problem setting, it follows from \eqref{u0pro} that
\[
u_0(t,x,i) = e^{q_{ii}(T-t)} w_0(t,x,i) = e^{q_{ii}(T-t)} V(x,r_i,\sigma_i,T-t,\alpha_i).
\]
We only derive the case $m=1$, as the further iterations with $m\in \{2,\cdots\}$ can be derived along the same line as \eqref{wmpro} and \eqref{umpro}.
With $m=1$, it holds by the Fubini theorem that
\begin{align*}
w_1(t,x,i)&= u_0(t,x,i) +\sum_{j\in\mathcal{M}\setminus\{i\}} q_{ij} \int_t^T e^{q_{ii}(s-t)} \mathbb{E}_0^{t,x,i}\left[ e^{-r_i(s-t)} w_0(s,X_s,j)\right] ds \\
&=u_0(t,x,i) +\sum_{j\in\mathcal{M}\setminus\{i\}} q_{ij} \int_t^T e^{q_{ii}(s-t)} \mathbb{E}_0^{t,x,i}\left[ e^{-r_i(s-t)} \mathbb{E}_0^{s,X_{s},j}\left[e^{-r_j(T-s)}g(X_T)\right]\right] ds,
\end{align*}
due to the non-negativity of the integrands.
Next, for each $s\in [0,T]$, we have %$X_T  \stackrel{\mathcal{L}}{=} X_s \exp[(r_j-\alpha_j-\sigma_j^2/2)(T-s)+\sigma_j\sqrt{T-s} Z]$ under $\mathbb{P}_0^{s,X_s,j}$, which yields
\[
\mathbb{E}_0^{s,X_{s},j}\left[e^{-r_j(T-s)}g(X_T)\right] = \mathbb{E}\left[e^{-r_j(T-s)}g\left( X_s \exp\left[(r_j-\alpha_j-\sigma_j^2/2)(T-s)+\sigma_j(W_T-W_s)\right]\right) \Big|\,\mathcal{F}_s \right],
\] 
%In a similar manner, it holds that for $s\in [t,T]$, %we have $X_s  \stackrel{\mathcal{L}}{=} x \exp [(r_i-\alpha_i-\sigma_i^2/2)(s-t)+\sigma_i\sqrt{s-t} Z]$ under $\mathbb{P}_0^{t,x,i}$.
%The iterative expectation above can be rewritten as
which yields
\begin{align*}
&\mathbb{E}_0^{t,x,i}\left[ e^{-r_i(s-t)} \mathbb{E}_0^{s,X_{s},j}\left[e^{-r_j(T-s)}g(X_T)\right]\right]\\
&\quad =\mathbb{E}_0^{t,x,i}\left[ e^{-r_i(s-t)} \mathbb{E}\left[e^{-r_j(T-s)}g\left( X_s \exp\left[(r_j-\alpha_j-\sigma_j^2/2)(T-s)+\sigma_j(W_T-W_s)\right]\right) \Big| \mathcal{F}_s \right]\right]\\
&\quad =\mathbb{E}\left[ e^{-r_i(s-t)} \mathbb{E}\left[e^{-r_j(T-s)}g\left(  x \exp\left[(r_i-\alpha_i-\sigma_i^2/2)(s-t)+\sigma_i(W_s-W_t)\right] \exp\left[(r_j-\alpha_j-\sigma_j^2/2)(T-s)+\sigma_j(W_T-W_s)\right]\right) \Big|\mathcal{F}_s \right]\right]\\
%&\quad =\mathbb{E}\left[ e^{-r_i(s-t) - r_j(T-s)} g\left(  x \exp\left[(r_i(s-t)+r_j(T-s)) - (\alpha_i(s-t)+\alpha_i(T-s)) - (\sigma_i^2(s-t) +\sigma_j^2(T-s))/2 + \sigma_i\sqrt{s-t} Z_1+\sigma_j\sqrt{T-s}Z_2\right]\right)\right]\\
&\quad =\mathbb{E}\left[ e^{-r(s)(T-t)} g\left(  x \exp\left[ (r(s)-\alpha(s) - \sigma^2(s)/2)(T-t) + \sigma(s)\sqrt{T-t} Z\right]\right)\right]=V(x,r(s),\sigma(s),T-t,\alpha(s)),
\end{align*}
where the third equality holds because the sum of two independent centered normal random variables with  variances $\sigma^2_i(s-t)$ and $\sigma^2_j(T-s)$ is identical in law to a centered normal random variable variance $\sigma^2_i(s-t)+\sigma^2_j(T-s) = \sigma^2(s) (T-t)$.
Note that the expectation from the third line on can be understood to be with respect to the underlying Brownian motion alone.
Therefore, we obtain
\begin{align*}
w_1(t,x,i) &= u_0(t,x,i) +\sum_{j\in\mathcal{M}\setminus\{i\}} q_{ij} \int_t^T e^{q_{ii}(s-t)} V(x,r(s),\sigma(s),T-t,\alpha(s)) ds,\\
u_1(t,x,i) &= u_0(t,x,i) +\sum_{j\in\mathcal{M}\setminus\{i\}} q_{ij} \int_t^T e^{q_{ii}(s-t)+q_{jj}(T-s)} V(x,r(s),\sigma(s),T-t,\alpha(s)) ds.
\end{align*}
The second and further iterations can be derived in a similar manner.

%\section*{Data availability statements}

%Data will be made available on reasonable request.

\end{document}